# A Complete Proof of the Strong Conjecture about $F$-Irregular Graphs

**Tatiana Dovzhenok**

Research Laboratory "Mathematics of Hybrid Intelligence Systems",
Francisk Skorina Gomel State University, Gomel, 246028, Belarus

and

**Artem Filuta**

Faculty of Applied Mathematics and Computer Science,
Belarusian State University, Minsk, 220030, Belarus

**Abstract.** A graph $G$ is called $F$-irregular if all its vertices have distinct $F$-degrees, defined as the number of subgraphs of $G$ isomorphic to a given graph $F$ and containing the respective vertex. We prove the Strong Conjecture about $F$-irregular graphs (Dovzhenok, Filuta, and Chuhai, 2024), which states that for every connected graph $F$ of order at least three, there exist infinitely many $F$-irregular graphs. Fundamentally generalizing the classical existence conjecture by Chartrand et al. (1987), this work presents an authorized English translation of our original February 2024 manuscript, which was publicly presented in full at two scientific conferences the same year.

# 1. Introduction

## 1.1. Research Priority and the Strong Conjecture

For graphs $F$ and $G$, the $F$*-degree* of a vertex $v$ in $G$ is the number of subgraphs of $G$ isomorphic to $F$ and containing $v$. The graph $G$ is called $F$-irregular if no two distinct vertices of $G$ have identical $F$-degrees. In the seminal paper [1] by Chartrand, Holbert, Oellermann, and Swart, the following conjecture was formulated:

**Conjecture 1** (Chartrand et al., 1987). *For every connected graph F of order three or more, there exists a non-trivial F-irregular graph.*

As follows from the 2021 monograph *Irregularity in Graphs* by A. Ali, G. Chartrand, and P. Zhang [2], Conjecture 1 had not been verified for any new infinite classes of graphs beyond the complete graphs and stars originally studied in [1].

In this paper, the authors present their independent results, obtained in February 2024.

To anchor our scientific priority, the 2024 Russian manuscript was deposited in open access on September 12, 2026 [3].

In that work, we proved a far more general statement than Conjecture 1, originally designated as the Generalized Conjecture and subsequently renamed by Dovzhenok, Filuta, and Chuhai [4] as the *Strong Conjecture about F-irregular graphs*:

**Conjecture 2** (The Strong Conjecture; Dovzhenok et al., 2024). *For every connected graph F of order three or more, there exist infinitely many F-irregular graphs.*

Section 2 of the present work contains an authorized translation of the manuscript [3], incorporating only minor technical corrections of typographical errors with no impact on the course of the proofs. This manuscript is the original competition entry (including its title page, introduction, and bibliography of that period), which accounts for the exhaustive level of detail in the presented proofs. Notably, all constructions and methods are presented in an accessible and straightforward manner.

## 1.2. Authors' Contributions and Approbations of Results

The foundational ideas and key methods of proof for simple graphs belong to Dovzhenok, while the verification of the constructions was performed jointly with Filuta (who was an 11th-grade student at that time). Under the scientific supervision of Dovzhenok, Filuta also obtained a proof of the statement analogous to Conjecture 2 for the class of multigraphs (Chapter 5 of Section 2). The results of the study were recorded by the authors in the manuscript [3], which was successfully presented in its entirety by Filuta (with Dovzhenok as the scientific supervisor) at the XXVIII Republican Contest of Scientific Research Papers of Students (Minsk, February 2024, First-Degree Diploma) and the XX Baltic Science and Engineering Fair (St. Petersburg, April 2024, Third-Degree Diploma).

The efficacy of the developed approaches has been confirmed by a series of peer-reviewed publications, yielding the following results:

- for each 2-connected graph $F$ with minimum degree $\delta(F) = 2$, an infinite family of $F$-irregular graphs was constructed jointly with Chuhai [4];
- regarding Conjecture 1, Dovzhenok completely settled the case of paths $P_n$ of arbitrary order $n \geq 3$ and constructed infinite families of $P_4$-irregular graphs [5];
- Dovzhenok also confirmed Conjecture 2 for an arbitrary graph $F$ of diameter 2 [6].

Furthermore, the ideas and constructions presented in [3] find their reflection in recent works by one of the authors [7], written jointly with Lukashenko and Filiuta, and [8], where the focus shifts toward the cases of oriented graphs and rooted $F$-irregularity, respectively

## 1.3. Technical Notes and Notation

In the present work, only finite and undirected graphs are considered. Simple graphs are defined as graphs without loops or multiple edges; multigraphs are defined as graphs that admit both loops and multiple edges. An edge with vertices $u$ and $v$ is denoted throughout the text in parentheses (u, v). The $F$-degree of a vertex $v$ in a graph $G$ is denoted by $Fdeg_G v$.

## 2. Authorized Translation of the 2024 Manuscript [3]

---

# Proof of the Conjecture About *F*-Irregular Graphs

**Author:**
Artem V. Filuta
Grade 11, State Educational Institution
"Secondary School No. 30 of Gomel"

**Research Advisor:**
Tatiana S. Dovzhenok,
PhD in Physics and Mathematics
Mathematics Teacher,
Secondary School No. 30 of Gomel

---

## 2.0. Introduction (2024)

$F$-irregular graphs are the object of study of this paper. The definition of an $F$-irregular graph was first introduced in [1] in 1987. In that paper, the Authors generalized the well-known graph-theoretic notions of vertex degree and irregular graph. We now recall these generalizations in more detail.

Consider two graphs $F$ and $G$ without loops or multiple edges. Let $v \in V(G)$.

Define the $F$-degree of a vertex $v$ in a graph $G$ as the number of subgraphs of $G$ that are isomorphic to $F$ and contain the vertex $v$.

We say that a graph $G$ is $F$-irregular if the $F$-degrees of its vertices are pairwise distinct.

In [1], the authors proved that for every natural number $n \geq 3$ there exist a $K_n$-irregular graph and a $K_{1,n-1}$-irregular graph, and posed the following **Conjecture**:

> *for every connected graph $F$ of order three or more, there exists a nontrivial $F$-irregular graph.*

Until 2024, this Conjecture remained an open problem in graph theory, although starting in 2021 the topic began to develop intensively, accumulating new results. Thus, in 2022 the existence of infinitely many $F$-irregular graphs was proved for $F \in \{K_n, K_{1,n-1}, C_n | n \geq 3\}$ [2]. The most interesting results, however, came in 2023: it was proved that the number of $F$-irregular graphs is infinite for any $F$ of order 3 or more containing a pendant vertex [3], and in [4] a result was published for graphs $F$ without cut vertices and with minimum vertex degree equal to 2. In addition, we formulated the **Generalized Conjecture about *F*-irregular graphs**:

> *for every connected graph $F$ of order three or more, there exist infinitely many $F$-irregular graphs.*

The initial goal of our study was to prove the Generalized Conjecture about $F$-irregular graphs. This proof is contained in Chapters 1–4 of the present paper. However, in the course of studying this problem we dealt only with simple graphs. This circumstance led us to an interesting question:

"Would a similar conjecture hold for graphs with loops and multiple edges, that is, for multigraphs?"

The answer to this question was successfully found and is given in Chapter 5.

The proof of the Conjecture consists of three stages:

1) The proof for complete graphs ($K_n$-Theorem);

2) The proof for graphs containing a pendant vertex (Pendant Theorem);

3) The proof for graphs of diameter 2 or more without pendant vertices (Fat Theorem).

As noted above, the $K_n$-Theorem was proved in [2]; however, in Chapter 1 we give a simplified, constructive proof of this theorem.

In Chapter 2, for the sake of completeness, we include the Pendant Theorem, previously proved in [3].

Chapter 3 contains a new result — the Fat Theorem (fat in a good sense).

In Chapter 4, building on the first three chapters, we formulate and prove the Simple Theorem, which, in essence, coincides with the Generalized Conjecture about $F$-irregular graphs.

Chapter 5 investigates graphs containing loops and multiple edges. The outcome of this investigation is the Multi-Theorem, which is the analogue of the Simple Theorem for multigraphs.

The proofs of all lemmas are given in the Appendix.

## 2. 1. Chapter 1. The $K_n$-Theorem (2024)

Let $n \in N, n \geq 3$. In this Chapter, the question of the cardinality of $K_n$-irregular graphs will be investigated.

### 2.1.1. The graph $A_{2l-1}$

Let $l \in N, l > 4$. Consider the graph $\boldsymbol{A_{2l-1}}$ with vertex set $V(A_{2l-1}) = \{1, 2, 3, \ldots, 2l - 1\}$, in which the vertices 1, 2, 3, …, $l$ form a complete graph, and the vertices $l + 1$, $l + 2$, …, $2l - 1$ also form a complete graph. Moreover, for every integer $i$ such that $l + 1 \leq i \leq 2l - 1$, the vertex $i$ is adjacent to the vertices $i - l + 1$, $i - l + 2$, …, $l$. There are no other edges in the graph $A_{2l-1}$.

For clarity, we depict the graph $A_{2l-1}$ as two levels: the upper level contains the vertices 1, 2, 3, …, $l$, and the lower level contains $l + 1$, $l + 2$, …, $2l - 1$, as shown in the figure. Then any two vertices within either level are adjacent, and every vertex of the lower level is adjacent to all vertices of the upper level lying strictly above it or to its right.

**Graph $\boldsymbol{A_{2l-1}}$**

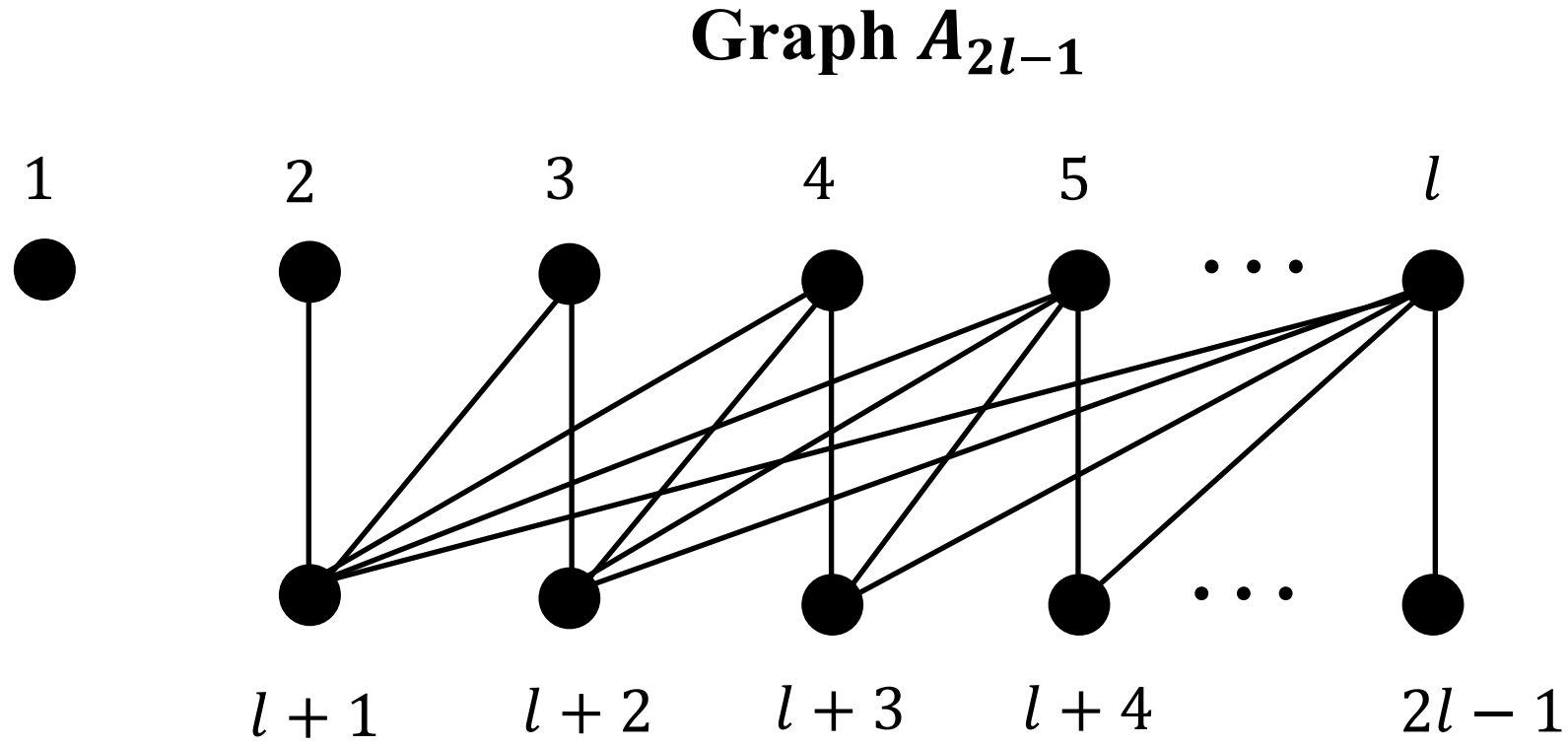

Let $i \in V(A_{2l-1})$. Denote by $a_i$ the $K_n$-degree of the vertex $i$ in the graph $A_{2l-1}$.

**Lemma 1.1.** *For every* $i \in \{1, 2, 3, \ldots, l - 1\}$ *the following equality holds:*

$$a_i = a_{2l-i}$$

**Lemma 1.2.** *Let $l > n$. For every $i \in \{1, 2, 3, \dots, l\}$ the following equality holds:*

$$a_i = \binom{l-1}{n-1} + (i-1)\binom{l-2}{n-2}.$$

## 2.1.2. The graph $B_{2l}$

Let $l \in N, l > 4$. Consider the graph $A_{2l-1}$ and add to it a vertex $2l$ adjacent to the vertices 1, 2, …, $l$. Call the resulting graph $\boldsymbol{B_{2l}}$.

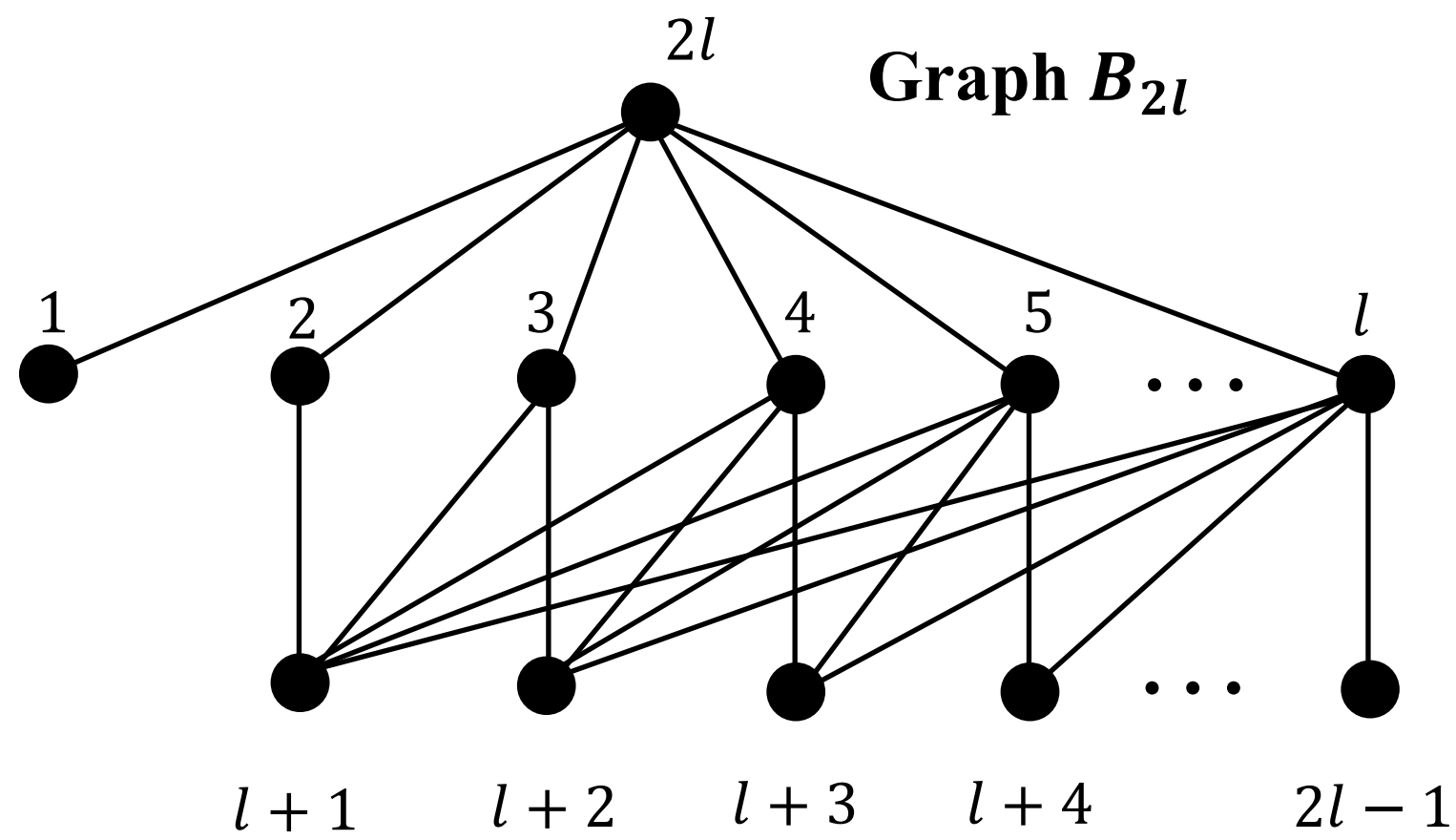

Let $i \in V(B_{2l})$. Denote by $b_i$ the $K_n$-degree of the vertex $i$ in the graph $B_{2l}$.

**Lemma 1.3.** *Let $l > n$. The $K_n$-degrees of the vertices of the graph $B_{2l}$ satisfy the following equalities:*

1) $b_{2l} = \binom{l}{n-1}$;

2) $b_i = \binom{l}{n-1} + (i-1)\binom{l-2}{n-2}$, *where* $i \in \{1, 2, 3, \dots, l\}$;

3) $b_i = \binom{l-1}{n-1} + (2l-i-1)\binom{l-2}{n-2}$, *where* $i \in \{l+1, l+2, \dots, 2l-1\}$.

## 2.1.3. The graph $U_{2l}$

Let $l \in N, l > n$. Consider the graph $B_{2l}$ and add to it the edge $(2l, 2l-n+2)$. Call the resulting graph $\boldsymbol{U_{2l}}$.

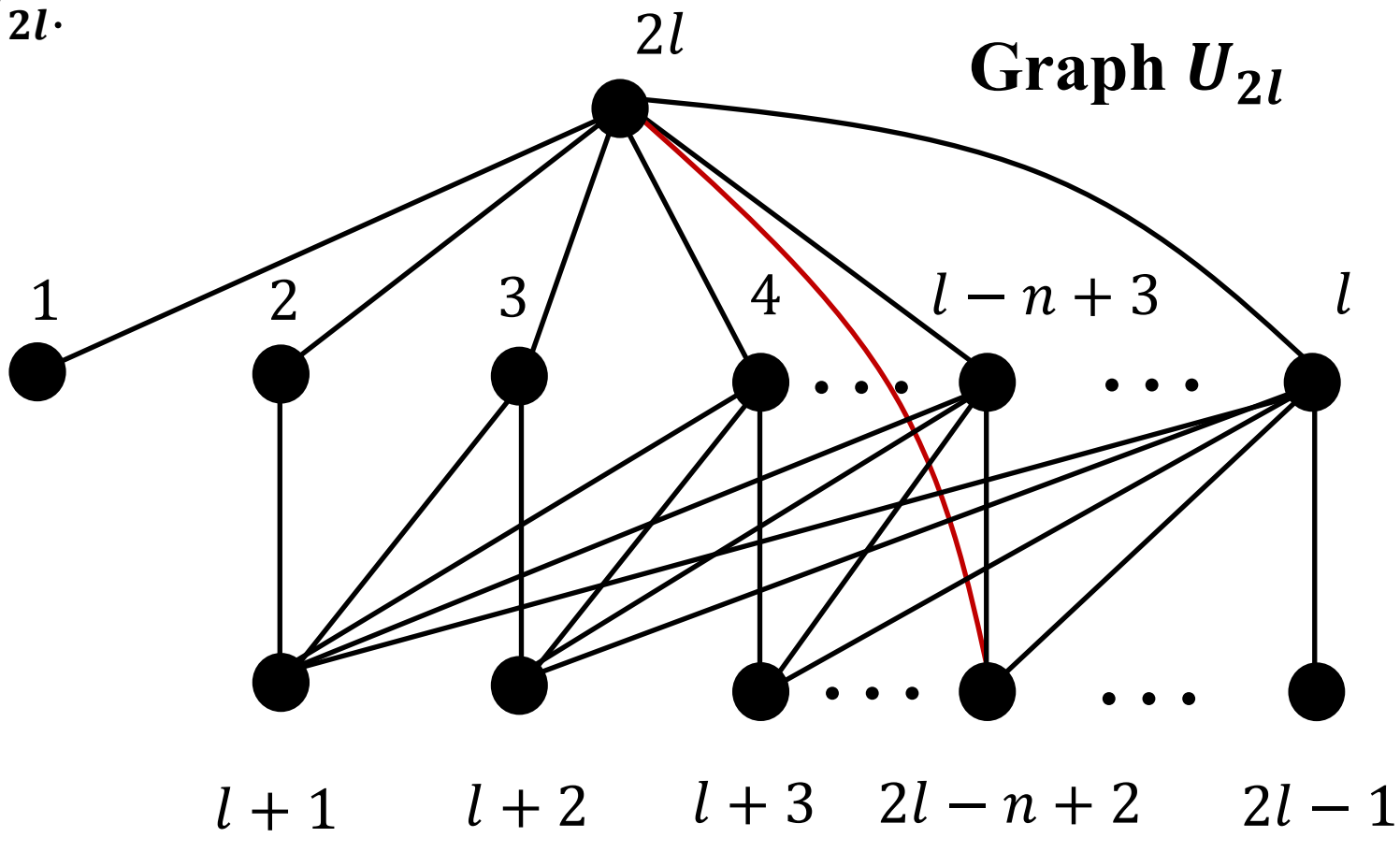

Let $i \in V(U_{2l})$. Denote by $u_i$ the $K_n$-degree of the vertex $i$ in the graph $U_{2l}$.

**Lemma 1.4**. *Let $l > n$. The $K_n$-degrees of the vertices of the graph $U_{2l}$ satisfy the following equalities:*

1) $u_{2l} = \binom{l}{n-1} + 1$;

2) $u_i = \binom{l}{n-1} + (i-1)\binom{l-2}{n-2},\ \ where\ \ i \in \{1, 2, 3, \dots, l-n+2\}$;

3) $u_i = \binom{l-1}{n-1} + (2l-i-1)\binom{l-2}{n-2},\ \ where\ \ i \in \{l+1, l+2, \dots, 2l-1\} \backslash \{2l-n+2\}$;

4) $u_i = \binom{l}{n-1} + (i-1)\binom{l-2}{n-2} + 1,\ \ where\ \ i \in \{l-n+3, l-n+4, \dots, l\}$;

5) $u_{2l-n+2} = \binom{l-1}{n-1} + (n-3)\binom{l-2}{n-2} + 1$.

**Lemma 1.5**. *Let $l > (n-2)! + 6,\ \ i \in \{1, 2, 3, \dots, l, 2l\},\ \ j \in \{l+1, l+2, \dots, 2l-1\}$. Then*

$$u_i \neq u_j$$

## 2.1.4. The $K_n$-Theorem

***$K_n$*-Theorem**. *For every natural number n ≥ 3, there exist infinitely many $K_n$-irregular graphs.*

**Proof of the $K_n$-Theorem.**

Let $l \in N, l > (n-2)! + 6$. Consider the graph $U_{2l}$. We prove that it is a $K_n$-irregular graph.

From Lemma 1.4 it follows that

$$u_1 < u_{2l} < u_2 < u_3 < \dots < u_l \quad \text{and} \quad u_{l+1} > u_{l+2} > \dots > u_{2l-1}.$$

Furthermore, from Lemma 1.5 we have

$$u_i \neq u_j \ \ for\ \ i \in \{1, 2, 3, \dots, l, 2l\},\ j \in \{l+1, l+2, \dots, 2l-1\}.$$

From the above we conclude that

$$u_a \neq u_b,\ \ for\ \ a, b \in \{1, 2, \dots, 2l\},\ a \neq b,$$

which precisely implies that $U_{2l}$ is a $K_n$-irregular graph.

Since there are infinitely many natural numbers $l > (n-2)! + 6$, it follows that there are also infinitely many $K_n$-irregular graphs. □

# 2.2. Chapter 2. The Pendant Theorem (2024)

Let $n \in N, n \geq 3$. In this Chapter we investigate the existence of $F$-irregular graphs in the case where $F$ is a connected graph of order $n$ containing at least one pendant vertex.

## 2.2.1. The graph $X_{2l}$

Let $l \in N, l > n$. Consider the graph $\boldsymbol{X_{2l}}$ with vertex set $V(X_{2l}) = \{1, 2, 3, \dots, 2l\}$, in which the vertices 1, 2, 3, …, $l$ form a complete graph, and the vertex $2l$ is adjacent to the vertex 1. Moreover, for every integer $i$ such that $l + 1 \leq i \leq 2l - 1$, the vertex $i$ is adjacent to the vertices $i - l + 1$, $i - l + 2$, …, $l$. There are no other edges in the graph $X_{2l}$.

For clarity, we depict the vertices $\{1, 2, 3, \dots, 2l-1\}$ of the graph $X_{2l}$ as two levels: the upper level contains the vertices 1, 2, 3, …, $l$, and the lower level contains $l+1$, $l+2$, …, $2l-1$, as shown in the figure. Then any two vertices in the upper level are adjacent, no two vertices in the lower level are adjacent, and every vertex of the lower level is adjacent to all vertices of the upper level lying strictly above it or to its right.

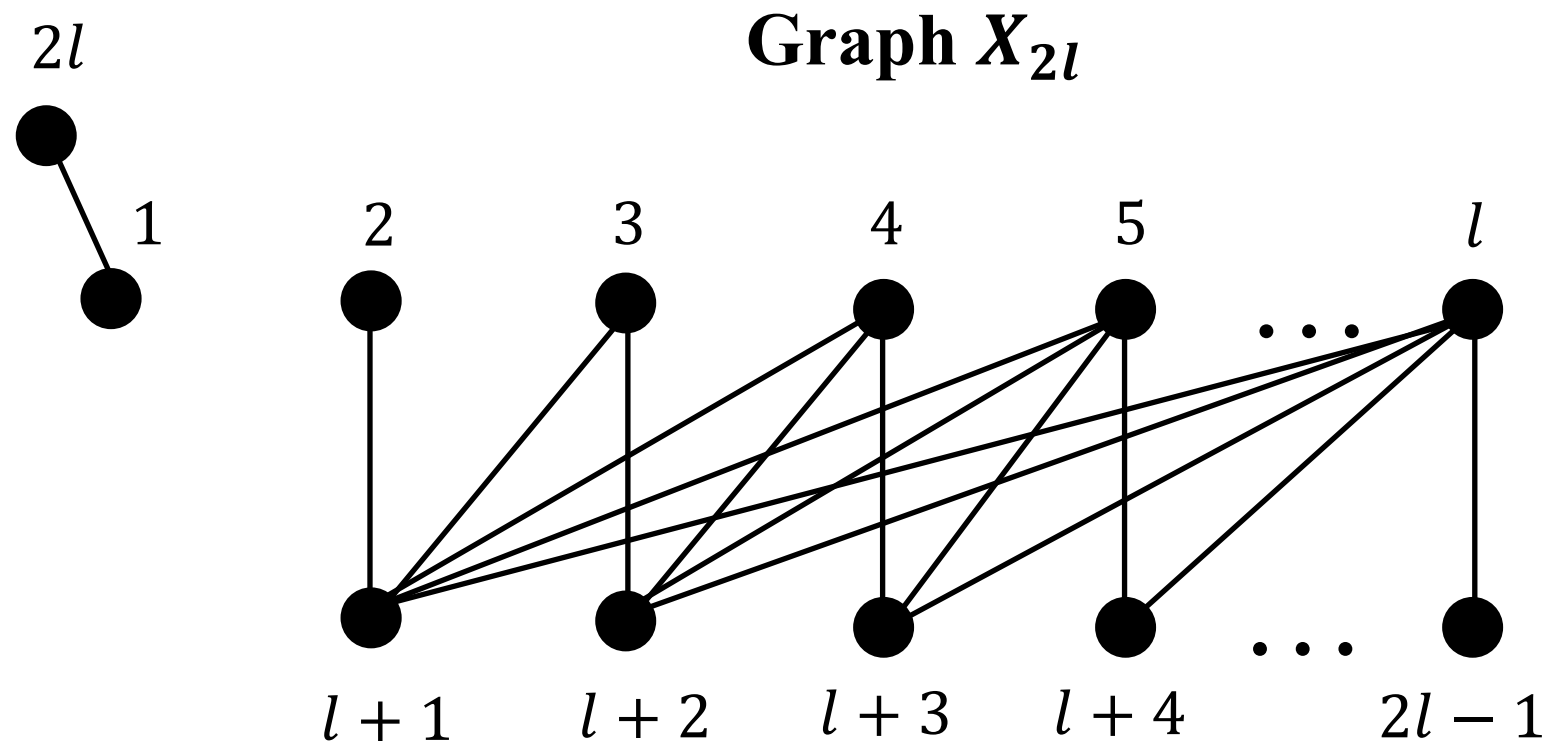

Let $i \in V(X_{2l})$. Denote by $x_i$ the $F$-degree of the vertex $i$.

**Lemma 2.1**. *Let $l > n$. For every $i \in \{2, 3, \dots, l-1\}$ the following inequality holds:*

$$x_{i+1} > x_i.$$

**Lemma 2.2**. *Let $l > n$. For every $i \in \{l+1, l+2, \dots, 2l-2\}$ the following inequality holds:*

$$x_i > x_{i+1}.$$

**Lemma 2.3.** *Let $l > n$. Then* $x_{l+1} < x_1$.

**Lemma 2.4**. *Let $l > n$, $\binom{l-2}{n-2} > n!\binom{2l-4}{n-3}$. Then $x_2 > x_1$.*

**Lemma 2.5**. *Let $> n$, $l \vdots 2$, and $\frac{l}{2} \cdot \binom{\frac{l}{2}-1}{n-3} > n!\binom{2l-4}{n-3}$. Then* $x_{2l-1} > x_{2l}$.

## 2.2.2. Two important inequalities

**Lemma 2.6**. *Let* $k_1, k_2, m_1, m_2, x \in \mathbb{N}$, $b_1, b_2 \in \mathbb{Z}$, $y \in \mathbb{Z}^+$, $x > y$, $r_1 > 0$, $r_2 \in \mathbb{R}$. *Then there exists a value $N_1$ such that for every natural number $l$ with $l \geq N_1$, $l \vdots \operatorname{lcm}(m_1, m_2)$, the following inequality holds:*

$$r_1 \cdot \binom{\frac{k_1}{m_1}l + b_1}{x} > r_2 \cdot \binom{\frac{k_2}{m_2}l + b_2}{y}.$$

**Lemma 2.7**. *Let* $k_1, k_2, m_1, m_2 \in \mathbb{N}$, $b_1, b_2 \in \mathbb{Z}$, $x \in \mathbb{Z}^+$, $p_1 > 0$, $r_1, r_2 \in \mathbb{R}$. *Then there exists a value $N_2$ such that for every natural number $l$ with $l \geq N_2$, $l \vdots \operatorname{lcm}(m_1, m_2)$, the following inequality holds:*

$$(p_1 l + r_1) \cdot \binom{\frac{k_1}{m_1}l + b_1}{x} > r_2 \cdot \binom{\frac{k_2}{m_2}l + b_2}{x}.$$

### 2.2.3. The Pendant Theorem

**Pendant Theorem**. *Let $n \in N, n \geq 3$ and let $F$ be a connected graph of order $n$ containing at least one pendant vertex. Then there exist infinitely many $F$-irregular graphs.*

**Proof of the Pendant Theorem.**
Consider the following **Condition**:

$$l \vdots 2, \quad l > n, \quad \binom{l-2}{n-2} > n!\binom{2l-4}{n-3}, \quad \frac{l}{2} \cdot \binom{\frac{l}{2}-1}{n-3} > n!\binom{2l-4}{n-3}.$$

We prove that infinitely many values of $l$ satisfy this **Condition**.

Indeed, by Lemmas 2.6 and 2.7, there exist numbers $N_1$ and $N_2$ such that for all natural numbers $l \geq N_1$ the inequality $\binom{l-2}{n-2} > n!\binom{2l-4}{n-3}$ holds, and for all even natural numbers $l \geq N_2$ the inequality $\frac{l}{2} \cdot \binom{\frac{l}{2}-1}{n-3} > n!\binom{2l-4}{n-3}$ holds. Therefore, the **Condition** holds for every even natural number $l \geq max(N_1, N_2, n+1)$. Clearly, there are infinitely many such $l$.

Consider the graph $X_{2l}$ for $l$ satisfying the above **Condition**.

From Lemmas 2.1–2.5 it follows that

$$x_l > x_{l-1} > \cdots > x_2 > x_1 > x_{l+1} > x_{l+2} > \cdots > x_{2l-1} > x_{2l}.$$

Thus, the graph $X_{2l}$ is $F$-irregular.

Since the **Condition** holds for infinitely many natural numbers $l$, there exist infinitely many graphs $X_{2l}$ that are $F$-irregular. □

## 2.3. Chapter 3. The Fat Theorem (2024)

Let $n \in \mathbb{N}$. In this Chapter we shall investigate the existence of $F$-irregular graphs in the case where $F$ is a connected graph of order $n$, with diameter at least 2 and no pendant vertices.

**Lemma 3.1.** *Let $n, t \in \mathbb{N}, n \geq 3$, let $F$ be a connected graph of order $n$ and diameter at least 2, and let $t$ be the minimum degree among the vertices of the graph $F$. Then there exist three vertices $u, v, w \in V(F)$ such that $degu = t$, $(u,v),(v,w) \in E(F)$, $(u,w) \notin E(F)$.*

### 2.3.1. The Graph $G_{2l-1}$

Let $l \in \mathbb{N}$, $l \vdots 2n^2$. Consider the graph $\boldsymbol{G_{2l-1}}$ with vertex set $V(G_{2l-1}) = \{1, 2, 3, \ldots, 2l-1\}$, in which the vertices $1, 2, 3, \ldots, l$ form a complete graph. Moreover, for every integer $i$ such that $l+1 \leq i \leq 2l-1$, the vertex $i$ is adjacent to the vertices $i-l+1, i-l+2, \ldots, l$.

In addition to these edges, the graph $G_{2l-1}$ contains a further set of edges, described below.

Let $t$ be the minimum degree among the vertices of the graph $F$, $t > 1$.

Consider a sequence of $\frac{t(t-1)}{2}$ vertices of the graph $G_{2l-1}$ labeled $1 + \frac{l}{n^2} \cdot j$, where $j \in \mathbb{Z}^+$, $j \leq \frac{t(t-1)}{2} - 1$. We call the vertices of this sequence red, and the remaining vertices of the graph $G_{2l-1}$, not belonging to this sequence, black. The existence of such a sequence follows from Lemma 3.2.

We join the first $t-1$ red vertices by edges to the vertex $2l-1$, the next $t-2$ vertices to the vertex $2l-2$, the next $t-3$ vertices to the vertex $2l-3$ and so on. Finally, we join the vertex labeled $1+\frac{l}{n^2}\cdot\left(\frac{t(t-1)}{2}-1\right)$ by an edge to the vertex $2l-t+1$. Then, in the graph $G_{2l-1}$, all vertices labeled $2l-t, 2l-t+1, \ldots, 2l-1$ have degree $t$, while the remaining vertices, by construction, have degrees greater than $t$. There are no other edges in the graph $G_{2l-1}$.

For clarity, we depict the vertices $\{1,2,3,\ldots,2l-1\}$ of the graph $G_{2l-1}$ as two levels: we place the vertices $1,2,3,\ldots,l$ on the upper level, and $l+1, l+2,\ldots,2l-1$ on the lower level, as shown in the figure. Then any two vertices on the upper level are adjacent, no two vertices on the lower level are adjacent, and every vertex on the lower level is adjacent to all vertices of the upper level located strictly above it or to its right. Moreover, the vertices $2l-t+1, 2l-t+2, \ldots, 2l-1$ are additionally joined by edges to the red vertices of the upper level in the manner described above.

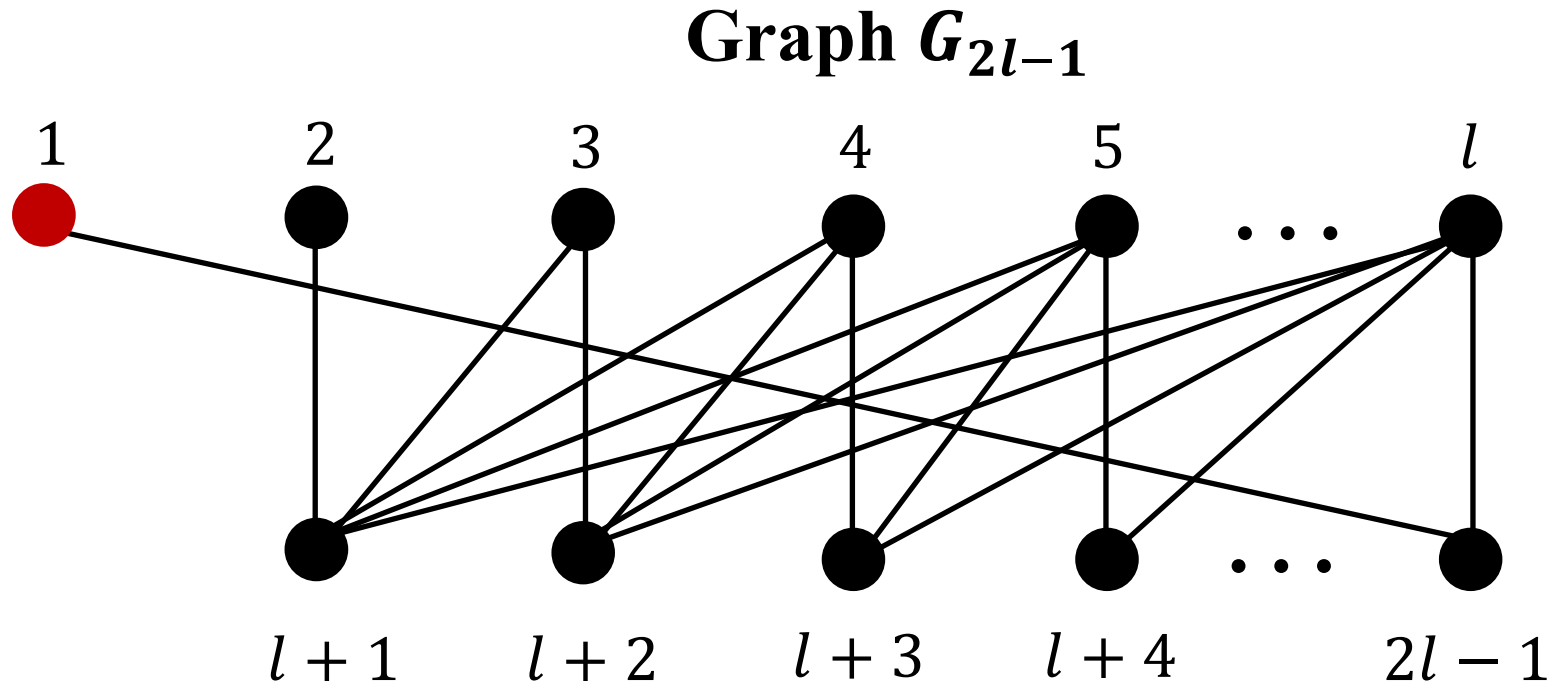

Let $i \in V(G_{2l-1})$. Denote by $g_i$ the $F$-degree of the vertex $i$.

**Lemma 3.2.** *Let $l \vdots 2n^2$. In the graph $G_{2l-1}$ the label of every red vertex does not exceed $\frac{l}{2}$.*

**Lemma 3.3.** *Let $l \vdots 2n^2$, $\binom{\frac{l}{2}-1}{n-2} > n!\binom{2l-t-3}{n-t-1}$. Then for every $i \in \{1,2,\ldots,l-1\}$ the following inequality holds:*

$$g_{i+1} > g_i.$$

**Lemma 3.4.** *Let $l \vdots 2n^2$. Then* $g_1 - g_{l+1} > 0$.

**Lemma 3.5.** *Let $l \vdots 2n^2$. For every $i \in \{l+1, l+2, \ldots, 2l-t-1\}$ the following inequality holds:*

$$g_i > g_{i+1}.$$

**Lemma 3.6.** *Let $l \vdots 2n^2$ and $\frac{l}{n^2}\binom{\frac{l}{n^2}}{n-t-2} > n!\,(t-1)\binom{2l-t-4}{n-t-2}$.*

*For every $i \in \{2l-t, 2l-t+1, \ldots, 2l-2\}$ the following inequality holds:*

$$g_i > g_{i+1}.$$

### 2.3.2. The Fat Theorem (fat in a good sense)

**Fat Theorem.** *Let* $n \in \mathbb{N}$, $n \geq 3$, $F$ *be* a *connected graph of order* n, *with diameter at least 2, having no pendant vertices. Then there exist infinitely many F-irregular graphs.*

**Proof of the Fat Theorem.**
Let $t$ be the minimum degree among the vertices of the graph $F$, $t > 1$. Then $t \leq n-1$.
Note that if $t = n-1$, then $F = K_n$ and, consequently, the diameter of the graph $F$ equals 1, which is false. Hence, $t \leq n-2$, whence $n-t-2 \geq 0$.
Moreover, $n-2 > n-t-1$ when $t > 1$.

Consider the following **Condition**:

$$l \in \mathbb{N},\ \ l \vdots 2n^2,\ \binom{\frac{l}{2}-1}{n-2} > n!\binom{2l-t-3}{n-t-1}\ ,\ \ \frac{l}{n^2}\binom{\frac{l}{n^2}}{n-t-2} > n!\,(t-1)\binom{2l-t-4}{n-t-2}.$$

We shall prove that this **Condition** is satisfied by infinitely many values of $l$.
Indeed, by Lemmas 2.6 and 2.7, there exist numbers $N_1$ and $N_2$ such that for every even natural number $l \geq N_1$ the inequality $\binom{\frac{l}{2}-1}{n-2} > n!\binom{2l-t-3}{n-t-1}$ holds, and for every natural number $l \geq N_2$ with $l \vdots n^2$ the inequality $\frac{l}{n^2}\binom{\frac{l}{n^2}}{n-t-2} > n!\,(t-1)\binom{2l-t-4}{n-t-2}$ *holds.*

Therefore, for every natural number $l \geq \max(N_1, N_2)$, $l \vdots 2n^2$ the **Condition** holds. Clearly, there are infinitely many such $l$.

Consider the graph $G_{2l-1}$ for $l$ satisfying the above **Condition**. By Lemmas 3.3–3.6, it follows that

$$g_l > g_{l-1} > \cdots > g_2 > g_1 > g_{l+1} > g_{l+2} > \cdots > g_{2l-1}.$$

Thus, the graph $G_{2l-1}$ is $F$-irregular. And since the **Condition** holds for infinitely many natural numbers $l$, there exist infinitely many graphs $G_{2l-1}$ that are $F$-irregular. ☐

## 2.4. Chapter 4. Proof of the Generalized Conjecture about *F*-irregular Graphs (2024)

**Simple Theorem.** *For every connected graph F of order 3 or more, there exist infinitely many F-irregular graphs***.**

**Proof of the Simple Theorem.**
Let $|F| = n \geq 3$ and $t$ be the minimum degree among the vertices of the graph $F$.

If $F = K_n$, then the result follows from the $K_n$-Theorem.

If $t = 1$, then the truth of the statement follows from the Pendant Theorem.

If $t > 1$ and $F \neq K_n$, then the graph $F$ has diameter 2 or more and contains no pendant vertices. Then the result follows from the Fat Theorem.

It remains to note that the three cases considered above completely cover all connected simple graphs of order 3 or more. ☐

## 2.5. Chapter 5. A Study of Graphs Containing Multiple Edges or Loops (2024)

In this Chapter we shall investigate the existence of $H$-irregular multigraphs in the case where $H$ is a connected multigraph.

**Definition 5.1.** *A **multigraph** is a graph containing at least one loop or at least one multiple edge.* (A simple graph does not fall under this definition.)

**Definition 5.2.** *The **skeleton** of a multigraph $H$ is a graph on the same vertex set as $H$, having no loops or multiple edges, in which two vertices $u$ and $v$ are adjacent if and only if they are adjacent in $H$.*

**Definition 5.3.** *Let $k \in \mathbb{N}$ and $H$ be a graph or multigraph in which the multiplicity of every loop and every edge does not exceed $k$. The **k-thickening** of a graph (multigraph) $H$* is *the multigraph, obtained from the skeleton H by adding $k$ loops to each vertex and increasing the multiplicity of each edge of the skeleton to $k$.*

**Multi-Theorem.** *For every connected multigraph $H$ there exist infinitely many nontrivial $H$-irregular multigraphs.*

**Proof of the Multi-Theorem.**
Suppose that $H$ is a multigraph with a single vertex having $a$ loops, $a \geq 1$, $a \in \mathbb{N}$. Then, an example of the desired multigraph is a multigraph $Z$, consisting of two vertices, 1 and 2, such that vertex 1 has $a$ loops and vertex 2 has $m > a$ loops. Since there are infinitely many natural values of $m > a$, there are infinitely many such desired multigraphs.

Now suppose that $|H| = 2$ and $V(H) = \{u, v\}$. Let the edge $(u, v)$ have multiplicity $t \in N$, vertex $u$ have $a$ loops, and vertex $v$ have $b$ loops, where $b \geq a$, $a + b + t \geq 2$. Then, an example of the desired multigraph is a graph $X$ on three vertices, 1, 2, 3, of the following form:

vertices 1 and 3 have $a$ loops, whereas vertex 2 has $c \geq b$ loops;
vertices 1, 2 are connected by an edge of multiplicity $t$, vertices 2, 3 are connected by an edge of multiplicity $t + 1$, and vertices 1, 3 are non-adjacent.

It is easy to see that

$$Hdeg_X 1 = \binom{c}{b}, \quad Hdeg_X 3 = \binom{c}{b}\binom{t+1}{t} = (t+1)\binom{c}{b},$$

$$Hdeg_X 2 = Hdeg_X 1 + Hdeg_X 3 = (t+2)\binom{c}{b}.$$

Therefore, $X$ is an $H$-irregular multigraph. Moreover, for different natural numbers $c \geq b$ we obtain different multigraphs.

Now suppose that $|H| \geq 3$. Let $t$ be the maximum edge multiplicity in the multigraph $H$, let p be the maximum number of loops attached to a single vertex, and let $k = \max(p, t)$.

Denote by $F$ the skeleton of the graph $H$. Since $H$ is connected, its skeleton $F$ is also connected. Therefore, by the Simple Theorem, there exist infinitely many $F$-irregular graphs.

We show how to transform each such graph $G$ into an $H$-irregular multigraph.

Let $D$ and $M$ be the $k$-thickenings of $H$ and $G$, respectively, and let $r$ be the number of subgraphs of the multigraph $D$ isomorphic to $H$. It is obvious that $r \geq 1$.

We prove that M is the desired $H$-irregular multigraph.

Let $u \in V(G)$, let $U$ be the set of all subgraphs of the graph $G$, isomorphic to $F$ and containing the vertex $u$, and let $U^+$ be the set of $k$-thickenings of all subgraphs from $U$.

Consider the map $f: U \to U^+$ which assigns to every graph $L \in U$ its $k$-thickening $L^+ \in U^+$. It is easy to see that the map $f$ is bijective. Therefore,

$$|U^+| = |U| = Fdeg_G u.$$

Let $u \in V(M) = V(G)$. Consider an arbitrary subgraph $X$ of $M$ that is isomorphic to $H$ and contains the vertex $u$. The skeleton $S$ of this multigraph $X$ is isomorphic to $F$ and contains the vertex $u$. Moreover $E(S) \subseteq E(G)$, by the definition of $M$ and $S$. Hence, $S \in U$ and, consequently, $X \subseteq S^+ \in U^+$. By the bijectivity of the map $f$, the subgraph $X$ belongs to exactly one $k$-thickening from $U^+$. Thus, the number of subgraphs from $M$ that are isomorphic to $H$ and containing the vertex $u$ is the sum of the numbers of such subgraphs in each of the $k$-thickenings from the set $U^+$. The required number of subgraphs in each of these $k$-thickenings equals $r$, since every $k$-thickening $L^+ \in U^+$ is isomorphic to the multigraph $D$ and, moreover, every subgraph from $L^+$ isomorphic to $H$ contains the vertex $u$. Thus, we have the equality

$$Hdeg_M u = r \cdot |U^+| = r \cdot Fdeg_G u.$$

If $u$ and $v$ are vertices of $G$, then

$$Fdeg_G u > Fdeg_G v \iff Fdeg_G u \cdot r > Fdeg_G v \cdot r \iff Hdeg_M u > Hdeg_M v.$$

Hence, if $G$ is $F$-irregular, then $M$ is an $H$-irregular multigraph. It remains to note that distinct graphs $G$ correspond to distinct multigraphs $M$. Therefore, there exist infinitely many $H$-irregular multigraphs. □

## 2.6. Conclusion (2024)

In this paper we have proved the **Generalized Conjecture about *F*-irregular graphs**:

*for every connected graph F of order 3 or more, there exist infinitely many F-irregular graphs*.

In particular, we have established the truth of the **Conjecture about *F*-irregular graphs [1]**:

*for every connected graph F of order 3 or more, there exists an F-irregular graph*.

Furthermore, we have formulated and proved the **Multi-Theorem**:

*for every connected multigraph H there exist infinitely many H-irregular multigraphs.*

The Multi-Theorem is an analog of the Simple Theorem about $F$-irregular graphs, but extended to multigraphs.

## 2.7. References (2024)

## 2.8. Appendix (2024)

**Proof of Lemma 1.1.** Let $i \in \{1, 2, 3, \dots, l-1\}$. We prove that in the graph $A_{2l-1}$, the vertices $i$ and $2l-i$ are symmetric. This can be easily established by placing the vertices of $A_{2l-1}$ on a straight line in ascending order of their labels from left to right. Then, the central vertex $l$ is adjacent to all other vertices. The vertex $i$ is adjacent to exactly $i-1$ vertices $1, 2, \dots, i-1$ to its left and to exactly $l-1$ vertices $i+1, i+2, \dots, l+i-1$ to its right, while the vertex $2l-i$ is adjacent to exactly $i-1$ vertices $2l-i+1, 2l-i+2, \dots, 2l-1$ to its right and to exactly $l-1$ vertices $l-i+1, l-i+2, \dots, 2l-i-1$ to its left. Thus, the vertices $i$ and $2l-i$ are symmetric in $A_{2l-1}$. Consequently, these vertices have equal $K_n$-degrees, meaning that $a_i = a_{2l-i}$ for any $i < l$. □

**Proof of Lemma 1.2**. Let us calculate $a_1$. To do this, we consider all subgraphs of $A_{2l-1}$ that are isomorphic to $K_n$ and contain the vertex 1. Since the vertex 1 is not adjacent to any vertex from the set $\{l+1, l+2, \dots, 2l-1\}$, such subgraphs can only include vertices from the set $\{1, 2, 3, \dots, l\}$. Since any two vertices from the set $\{1, 2, 3, \dots, l\}$ are adjacent in the graph $A_{2l-1}$, exactly one complete subgraph of $A_{2l-1}$ can be formed from the vertex 1 and any $n-1$ vertices of the set $\{2, 3, \dots, l\}$. Thus, $a_1 = \binom{l-1}{n-1}$.

Consider two consecutive vertices $i$ and $i+1$ in $A_{2l-1}$, where $i \in \{1, 2, 3, \dots, l-1\}$.

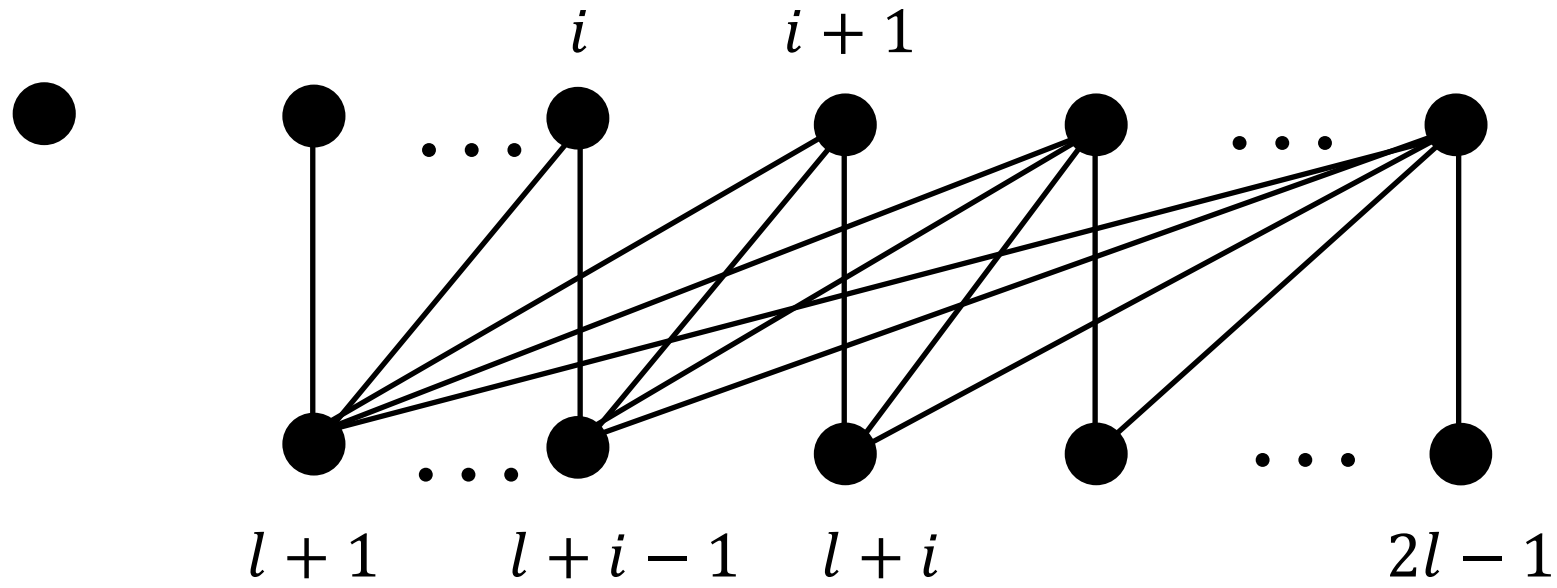

Let us find $a_{i+1} - a_i$. Notice that if we remove the edge $(i+1, l+i)$ from the graph $A_{2l-1}$, the vertices $i$ and $i+1$ become symmetric. This is because they are both adjacent to the vertices $l+1, l+2, \dots, l+i-1$, as well as to all vertices of the upper level other than $i$ and $i+1$. Furthermore, the vertices $i$ and $i+1$ are not adjacent to the vertices $l+i, \dots, 2l-1$ in the graph $A_{2l-1} \setminus (i+1, l+i)$. Consequently, the $K_n$-degrees of the vertices $i$ and $i+1$ in the graph $A_{2l-1} \setminus \{(i+1, l+i)\}$ coincide. Therefore, the difference between the $K_n$-degrees of these vertices in the graph $A_{2l-1}$ is equal to the number of subgraphs isomorphic to $K_n$ that contain the edge $(i+1, l+i)$ and do not contain the vertex $i$. Notice that each such subgraph contains no vertices from the set $\{l+i+1, l+i+2, \dots, 2l-1\}$, since none of these vertices is adjacent to the vertex $i+1$. Similarly, each such subgraph contains no vertices from the set $\{1, 2, \dots, i\}$, since none of these vertices is adjacent to the vertex $l+i$. Hence, the remaining vertices of such subgraphs must belong to the set $\{i+2, \dots, l+i-1\}$. And since any two vertices of the set $\{i+1, i+2, \dots, l+i\}$ are adjacent in the graph $A_{2l-1}$, the edge $(i+1, l+i)$ can be extended to a complete graph $K_n$ by choosing any $n-2$ vertices from the set $\{i+2, \dots, l+i-1\}$. This set of vertices can be selected in exactly $\binom{l-2}{n-2}$ ways. So, $a_{i+1} - a_i = \binom{l-2}{n-2}$. Given that $a_1 = \binom{l-1}{n-1}$, this yields:

$$a_i = \binom{l-1}{n-1} + (i-1)\binom{l-2}{n-2}, \; i \in \{1, 2, 3, \dots, l\}. \; \square$$

**Proof of Lemma 1.3.**
**3)** Let $i \in \{l+1, l+2, \dots, 2l-1\}$. Since the vertices $i$ and $2l$ are not adjacent in the graph $B_{2l}$, there exists no subgraph of $B_{2l}$ isomorphic to $K_n$ that contains both vertices $i$ and $2l$ simultaneously. Consequently, the $K_n$-degree of the vertex $i$ remains unchanged upon transitioning from the graph $A_{2l-1}$ to $B_{2l}$. Thus, by Lemmas 1.1 and 1.2, we obtain:

$$b_i = a_i = a_{2l-i} = \binom{l-1}{n-1} + (2l-i-1)\binom{l-2}{n-2} \ \text{ for } \ i \in \{l+1, l+2, \dots, 2l-1\}.$$

**1)–2)** Let $i \in \{1, 2, \dots, l\}$. Let us determine how the $K_n$-degree of the vertex $i$ changes upon transitioning from the graph $A_{2l-1}$ to the graph $B_{2l}$. To do this, we find the number of subgraphs isomorphic to $K_n$ that contain both vertices $i$ and $2l$ simultaneously. It is obvious that none of the vertices from the set $\{l+1, l+2, \dots, 2l-1\}$ belongs to such subgraphs. On the other hand, all vertices $1, 2, \dots, l, 2l$ form a complete subgraph of $B_{2l}$. Therefore, the number of the desired subgraphs is equal to the number of subsets of $n-2$ vertices from the set $\{1,2,3,\dots,l\}\backslash\{i\}$. There are exactly $\binom{l-1}{n-2}$ such subsets. Consequently, $b_i = a_i + \binom{l-1}{n-2}$. By Lemma 1.2, we obtain:

$$b_i = \binom{l-1}{n-1} + (i-1)\binom{l-2}{n-2} + \binom{l-1}{n-2} = \binom{l}{n-1} + (i-1)\binom{l-2}{n-2} \ \text{ for } i \in \{1, 2, \dots, l\}.$$

It remains to note that in the graph $B_{2l}$, the vertices 1 and $2l$ are symmetric. Therefore, $b_{2l} = b_1 = \binom{l}{n-1}$. □

**Proof of Lemma 1.4.** Consider the subgraphs isomorphic to $K_n$ that appeared in the graph $U_{2l}$ and were not present in the graph $B_{2l}$. It is obvious that all such subgraphs must contain the edge $(2l, 2l-n+2)$, and the remaining $n-2$ vertices must be adjacent to the vertices $2l$ and $2l-n+2$, as well as to each other. However, in the graph $U_{2l}$, only the $n-2$ vertices labeled $l-n+3, l-n+4, \dots, l$ satisfy this condition. Consequently, in $U_{2l}$, there exists exactly one subgraph isomorphic to $K_n$ that contains the edge $(2l, 2l-n+2)$. Therefore, the $K_n$-degrees of the vertices $l-n+3, l-n+4, \dots, l, 2l$, and $2l-n+2$ increase by 1 upon transitioning from $B_{2l}$ to $U_{2l}$, while the $K_n$-degrees of the remaining vertices remain unchanged. By Lemma 1.3, we obtain the required result. □

**Proof of Lemma 1.5.** Notice that the condition $l > (n-2)! + 6$ implies that $l > 2n$. Indeed, for $n = 3$, we get: $l > (3-2)! + 6 = 7 > 6$, and for $n \geq 4$, the following estimate holds:

$$l > (n-2)! + 6 \geq (n-2)(n-3) + 6 \geq 2(n-3) + 6 = 2n.$$

Let $i \in \{1, 2, 3, \dots, l, 2l\}$ and $j \in \{l+1, l+2, \dots, 2l-1\}$. We prove that $u_i \neq u_j$.

Assume the contrary, that is, $u_i = u_j$. It follows from Lemma 1.4 that

$$u_i = C_l^{n-1} + (p-1)C_{l-2}^{n-2} + m_p, \quad u_j = C_{l-1}^{n-1} + (2l-j-1)C_{l-2}^{n-2} + n_j,$$

where $m_p \in \{0, 1\}$, $p = i$ for $i \neq 2l$ and $p = 1$ for $i = 2l$; $n_j = 1$ for $j = 2l-n+2$ and $n_j = 0$ for $j \neq 2l-n+2$.

Then the equality $u_i = u_j$ is equivalent to the equality

$$\binom{l}{n-1} + (p-1)\binom{l-2}{n-2} + m_p = \binom{l-1}{n-1} + (2l-j-1)\binom{l-2}{n-2} + n_j \quad \Leftrightarrow$$

$$m_p + \binom{l-1}{n-2} = (2l-p-j)\binom{l-2}{n-2} + n_j \quad \Leftrightarrow$$

$$m_p + \frac{(l-1)!}{(n-2)!(l-n+1)!} = n_j + (2l-p-j)\frac{(l-2)!}{(n-2)!(l-n)!}.$$

Next, three cases are possible.

**1) $m_p = n_j$.** Then the equality under investigation takes the form:

$$\frac{(l-1)!}{(n-2)!(l-n+1)!} = (2l-p-j)\frac{(l-2)!}{(n-2)!(l-n)!} \Leftrightarrow \frac{l-1}{l-n+1} = 2l-p-j.$$

Since $(2l-p-j)$ is an integer, then $\frac{l-1}{l-n+1}$ is also an integer. Furthermore, $l-1 > l-n+1 > 0$ for $n \geq 3$ and $l > 2n$. Hence, $\frac{l-1}{l-n+1} \geq 2$, which implies that $2n-3 \geq l$. This contradicts the inequality $l > 2n$. Thus, Case 1 is impossible.

**2) $m_p = 1, n_j = 0$.** Then, the equality under investigation takes the form:

$$1 + \frac{(l-1)!}{(n-2)!(l-n+1)!} = (2l-p-j)\frac{(l-2)!}{(n-2)!(l-n)!} \quad \Leftrightarrow$$

$$(n-2)! + (l-n+2)(l-n+3)\cdot \ldots \cdot (l-1) =$$

$$= (2l-p-j)(l-n+1)(l-n+2)\cdot \ldots \cdot (l-2).$$

If $n = 3$, then $l = (2l-p-j)(l-2)$. Since $(2l-p-j)(l-2) \vdots l-2$, it follows that $l \vdots l-2$, whence $2 = l-(l-2) \vdots l-2$, which is false for $l > 2n = 6$.

If $n \geq 4$, then $(l-n+2)(l-n+3)\cdot \ldots \cdot (l-1) \vdots l-2$ and $(2l-p-j)(l-n+1)(l-n+2)\cdot \ldots \cdot (l-2) \vdots l-2$. Consequently, $(n-2)! \vdots l-2$. Therefore, $(n-2)! \geq l-2$, which contradicts the condition $l > (n-2)! + 6$.

**3) $m_p = 0, n_j = 1$.** Here $j = 2l-n+2$ and the equality under study will take the form:

$$\frac{(l-1)!}{(n-2)!(l-n+1)!} = 1 + (n-2-p)\frac{(l-2)!}{(n-2)!(l-n)!} \quad \Leftrightarrow$$

$$(l-n+2)(l-n+3)\cdot \ldots \cdot (l-1) = (n-2)! +$$

$$+(n-2-p)(l-n+1)(l-n+2)\cdot \ldots \cdot (l-2).$$

If $n = 3$, then $l-1 = 1 + (1-p)(l-2)$, whence $p = 0$. This is a contradiction.
If $n \geq 4$, then $(l-n+2)(l-n+3)\cdot \ldots \cdot (l-1) \vdots l-2$ and $(n-2-p)(l-n+1)(l-n+2)\cdot \ldots \cdot (l-2) \vdots l-2$. Hence, $(n-2)! \vdots l-2$.
So, $(n-2)! \geq l-2$, which contradicts the condition $l > (n-2)! + 6$.
Thus, in each of the cases 1), 2), and 3), a contradiction is obtained. Consequently, $u_i \neq u_j$ for any $i \in \{1,2,3,\ldots,l,2l\}$ and $j \in \{l+1, l+2, \ldots, 2l-1\}$. □

**Proof of Lemma 2.1**. Let $l > n$. Consider in $X_{2l}$ two consecutive vertices $i$ and $i+1$, where $i \in \{2,3,\ldots,l-1\}$. We estimate the difference $x_{i+1} - x_i$.

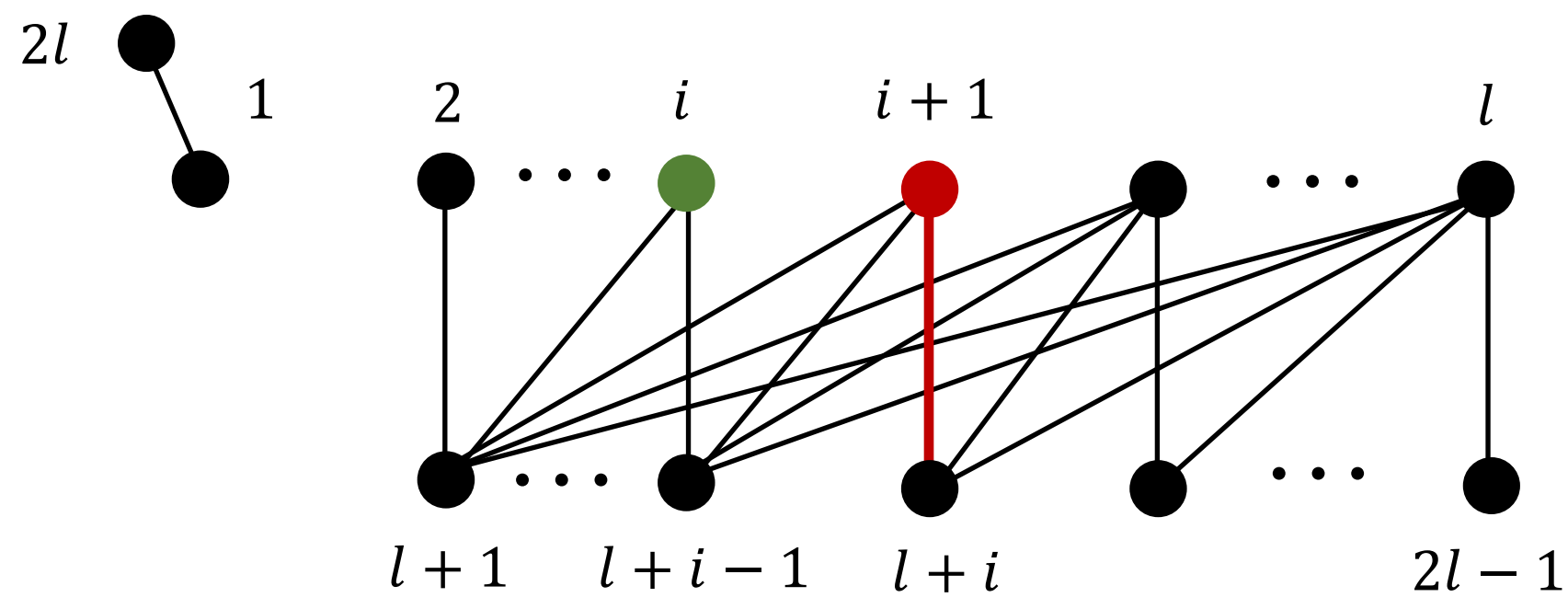

Note that in the graph obtained by deleting the edge $(i+1, l+i)$ from $X_{2l}$, the vertices $i$ and $i + 1$ are symmetric, since both are adjacent to the vertices $l+1, l+2, \ldots, l+i-1$, as well as to all vertices of the upper level other than $i$ and $i+1$. Moreover, the vertices $i$ and $i+1$ in the graph $X_{2l} \setminus \{(i+1, l+i)\}$ are not adjacent to the vertices $l+i, \ldots, 2l-1, 2l$. Consequently, the $F$-degrees of the vertices $i+1$ and $i$ in $X_{2l} \setminus (i+1, l+i)$ coincide. Hence the difference of the $F$-degrees $(x_{i+1} - x_i)$ of these vertices in $X_{2l}$ equals the number of subgraphs isomorphic to $F$ that contain the edge $(i+1, l+i)$ and do not contain the vertex $i$. The set of such subgraphs is nonempty. Let us prove this.

Consider an arbitrary set of $n-1$ vertices on the upper level containing the vertex $i+1$ and not containing the vertex $i$. Since any two vertices of such a set are adjacent in $X_{2l}$, from all these vertices, the vertex $l+i$, and the edge $(i+1, l+i)$, by adding the necessary edges one can form at least one subgraph isomorphic to $F$ and containing the edge $(i+1, l+i)$, in which the role of the pendant vertex is played by the vertex $l+i$, adjacent to the vertex $i+1$. By construction, such a subgraph of $X_{2l}$ does not contain the vertex $i$. Consequently, $x_{i+1} - x_i > 0$. □

**Proof of Lemma 2.2**. Let $l > n$. Consider in $X_{2l}$ two consecutive vertices $i$ and $i+1$, where $i \in \{l+1, l+2, \ldots, 2l-2\}$. We estimate the difference $x_i - x_{i+1}$. If we delete the edge $(i, i-l+1)$ from the graph $X_{2l}$, then the vertices $i$ and $i+1$ become symmetric, since both are adjacent to the vertices $i-l+2, i-l+3, \ldots, l$, and are adjacent to no other vertices. Consequently, the $F$-degrees of the vertices $i$ and $i+1$ in the graph $X_{2l} \setminus \{(i, i-l+1)\}$ coincide.

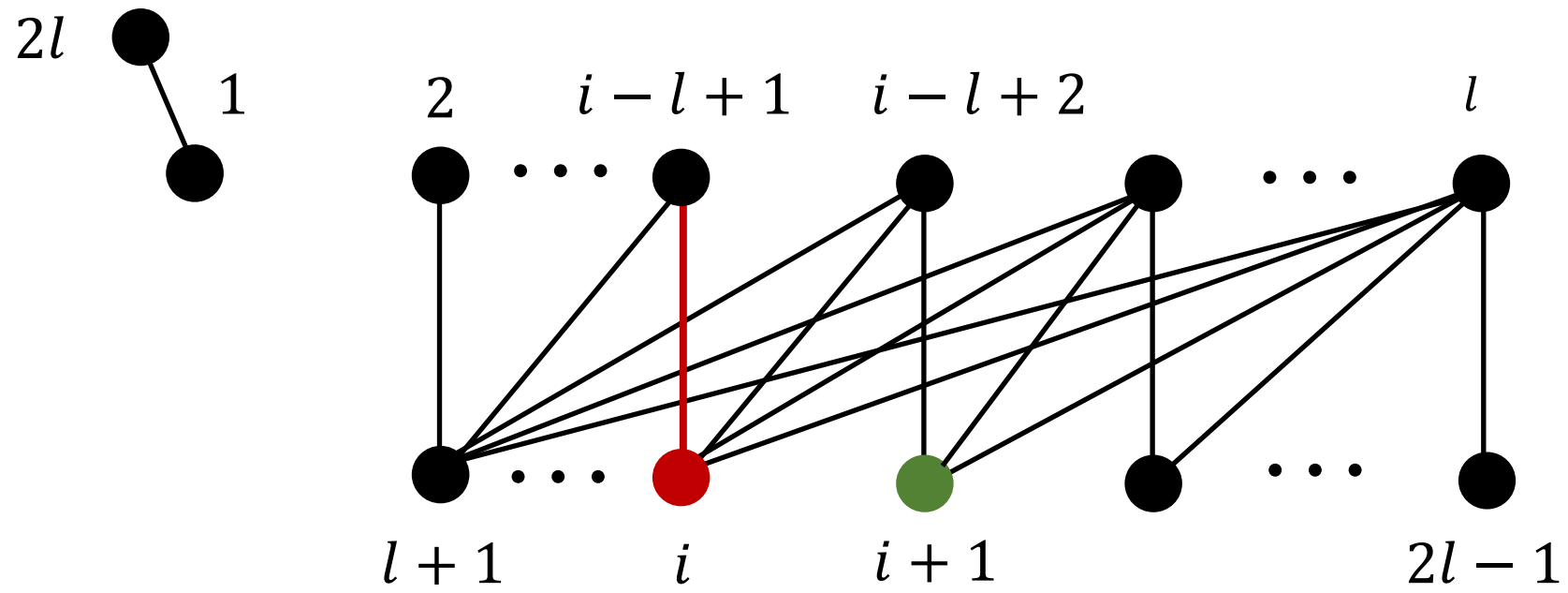

Hence the difference of the $F$-degrees $(x_i - x_{i+1})$ of these vertices in $X_{2l}$ equals the number of subgraphs isomorphic to $F$ that contain the edge $(i, i-l+1)$ and do not contain the vertex $i+1$. The set of such subgraphs is nonempty. Let us prove this.

Consider an arbitrary set of $n-1$ vertices of the upper level containing the vertex $i-l+1$. Since any two vertices of such a set are adjacent in the graph $X_{2l}$, from all these vertices, the vertex $i$, and the edge $(i, i-l+1)$, by adding the necessary edges one can form at least one subgraph isomorphic to $F$ and containing the edge $(i, i-l+1)$, in which the role of the pendant vertex is played by the vertex $i$, adjacent to the vertex $i-l+1$. By construction, such a subgraph of $X_{2l}$ does not contain the vertex $i+1$. Consequently, $x_i - x_{i+1} > 0$. □

**Proof of Lemma 2.3**. Let $l > n$. Consider in the graph $X_{2l}$ two vertices 1 and $l+1$. Let us estimate the difference $x_1 - x_{l+1}$. If we delete from the graph $X_{2l}$ the edge $(1, 2l)$, then the vertices 1 and $l+1$ become symmetric, since in the resulting graph both are adjacent to the vertices $2, 3, \ldots, l$, and are adjacent to no other vertices. Consequently, the $F$-degrees of the vertices 1 and $l+1$ in the graph $X_{2l} \setminus \{(1, 2l)\}$ coincide. Hence the difference of the $F$-degrees $(x_1 - x_{l+1})$ of these vertices in the graph $X_{2l}$ equals the number of subgraphs isomorphic to $F$ that contain the edge $(1, 2l)$ and do not contain the vertex $l+1$.

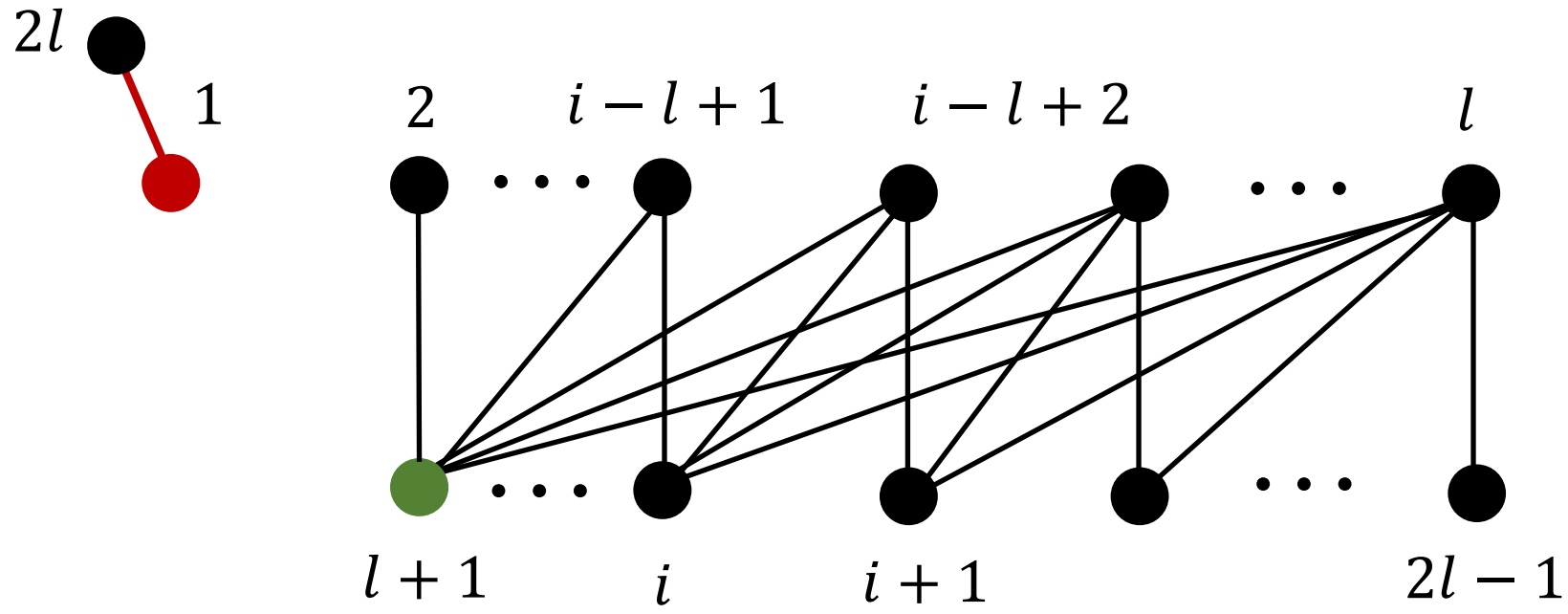

The set of such subgraphs is nonempty. Let us prove this.

Consider an arbitrary set of $n-1$ vertices on the upper level containing the vertex 1. Since any two vertices of such a set are adjacent in $X_{2l}$, from all these vertices, the vertex $2l$, and the edge $(1,2l)$, by adding the necessary edges one can form at least one subgraph isomorphic to $F$ and containing the edge $(1,2l)$, in which the role of the pendant vertex is played by the vertex $2l$, adjacent to the vertex 1. By construction, such a subgraph of $X_{2l}$ does not contain the vertex $l+1$. Consequently, $x_1 - x_{l+1} > 0$, whence $x_{l+1} < x_1$. □

**Proof of Lemma 2.4.** Let $l > n$ and $\binom{l-2}{n-2} > n!\binom{2l-4}{n-3}$. Consider in $X_{2l}$ two vertices 1 and 2. We prove that $x_2 > x_1$.

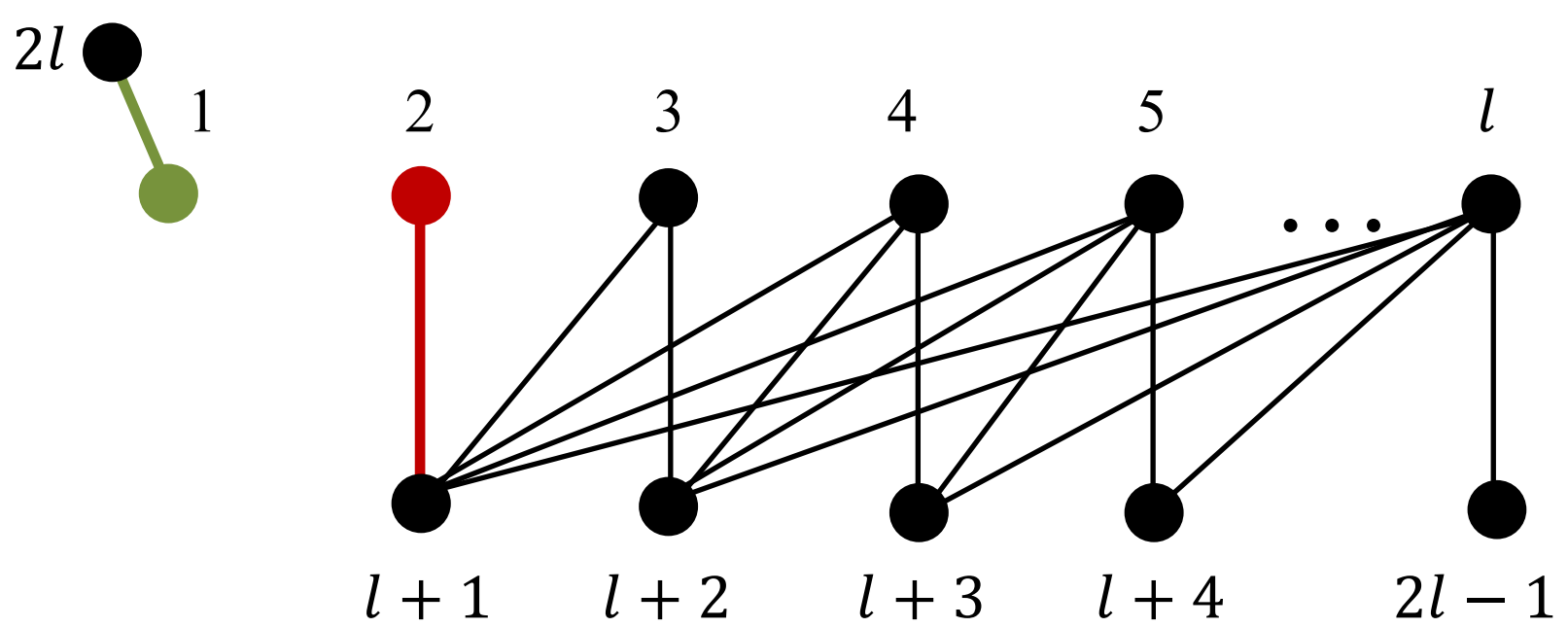

Let us estimate the difference $x_2 - x_1$. Upon deleting the edges $(1,2l)$ and $(2,l+1)$ from the graph $X_{2l}$, the vertices 1 and 2 become symmetric, since both are adjacent to the vertices 3, 4,…, $l$ and to each other, and are adjacent to no other vertices. Consequently, the $F$-degrees of the vertices 2 and 1 in the graph $X_{2l} \setminus \{(1,2l),\ (2,l+1)\}$ coincide. Hence the difference of the $F$-degrees $(x_2 - x_1)$ of these vertices in $X_{2l}$ equals the difference between the number of subgraphs isomorphic to $F$ that contain the edge $(2,l+1)$ and do not contain the vertex 1 (denote the set of such subgraphs by $M$), and the number of subgraphs isomorphic to $F$ that contain the edge $(1,2l)$ and do not contain the vertex 2 (denote the set of such subgraphs by $K$). Therefore, to prove the inequality $x_2 - x_1 > 0$ it suffices to show that $|M| > |K|$.

To this end, represent $K$ and $M$ as unions of disjoint sets:

$$K = K_1 \cup K_2 \quad and \quad M = M_1 \cup M_2$$

where

$$K_1 = \{G \in K \mid l+1 \notin V(G)\}, \quad K_2 = \{G \in K \mid l+1 \in V(G)\},$$

$$M_1 = \{H \in M \mid deg_H(2) = 1\}, \quad M_2 = \{H \in M \mid deg_H(2) > 1\},$$

and we prove that $|M_1| \geq |K_1|$ and $|M_2| > |K_2|$.

**1) We prove that $|M_1| \geq |K_1|$.**

Consider a mapping of the set $K_1$ into the set $M_1$, under which to each subgraph $G \in K_1$ with vertex set $V(G)$ and edge set $E(G)$ we assign a subgraph of $M_1$ according to the following rule:

$$2l \to 2, \quad 1 \to l+1, \quad i \to i \;\; for \;\; i \in V(G) \setminus \{1, 2l\},$$
$$(1, 2l) \to (l+1, 2), \quad (1, i) \in E(G) \to (l+1, i) \;\; for \;\; i \neq 2l,$$
$$(i, j) \in E(G) \to (i, j) \;\; for \;\; i, j \neq 1.$$

Clearly, this mapping is injective, so $|M_1| \geq |K_1|$.

**2) We prove that $|M_2| > |K_2|$.**

*Let us estimate $|K_2|$ from above.*

Every subgraph in $K_2$ contains the vertices 1, $l+1$, $2l$ and does not contain the vertex 2. The remaining $n-3$ vertices for such subgraphs can be chosen in $\binom{2l-4}{n-3}$ ways. Since to each set of $n$ vertices of the graph $X_{2l}$ there correspond at most $n!$ subgraphs isomorphic to $F$ on this vertex set, it follows that

$$|K_2| \leq n! \binom{2l-4}{n-3}.$$

*Let us estimate $|M_2|$ from below.*

Every subgraph in $M_2$ contains the vertices 2, $l+1$ together with the edge $(2, l+1)$ and does not contain the vertex 1. Next, consider only those subgraphs in $M_2$ for which the remaining $n-2 \geq 1$ vertices lie on the upper level of $X_{2l}$ and are distinct from 1, 2, and for which the vertex $l+1$ is a pendant vertex. Observe that, by the connectedness of the graph $F$, the vertex 2, adjacent to the vertex $l+1$, cannot be a pendant vertex.

Let $m$ be the number of such subgraphs. Note that since any two vertices on the upper level are adjacent in $X_{2l}$, by adding an arbitrary set of $n-2$ vertices on the upper level (distinct from 1 and 2) to the vertices 2, $l+1$ and the edge $(2, l+1)$, together with the necessary edges among the vertices of the chosen set, it is always possible to form at least one subgraph isomorphic to $F$ in which the vertex $l+1$ is pendant (adjacent to 2) and which does not contain the vertex 1. Since there are in total $\binom{l-2}{n-2}$ ways to choose $n-2$ vertices from the $l-2$ admissible vertices on the upper level, it follows that $m \geq \binom{l-2}{n-2}$. Consequently,

$$|M_2| \geq m \geq \binom{l-2}{n-2}.$$

Taking into account the assumption $\binom{l-2}{n-2} > n! \binom{2l-4}{n-3}$ and the inequality $|K_2| \leq n! \binom{2l-4}{n-3}$ proved above, we obtain $|M_2| > |K_2|$.

From 1), 2) it follows that $|M| > |K|$. □

**Proof of Lemma 2.5.** Let $l > n$, $l \,\vdots\, 2$, $\frac{l}{2}\binom{\frac{l}{2}-1}{n-3} > n! \binom{2l-4}{n-3}$. Consider in the graph $X_{2l}$ two vertices $2l-1$ and $2l$. We estimate the difference $x_{2l-1} - x_{2l}$. Clearly, it equals the difference between the number of subgraphs isomorphic to $F$, containing the vertex $2l-1$ and not containing the vertex $2l$ (we denote the set of such subgraphs by $A$), and the number of subgraphs isomorphic to $F$, containing the vertex $2l$ and not containing the vertex $2l-1$ (we denote the set of such subgraphs by $B$).

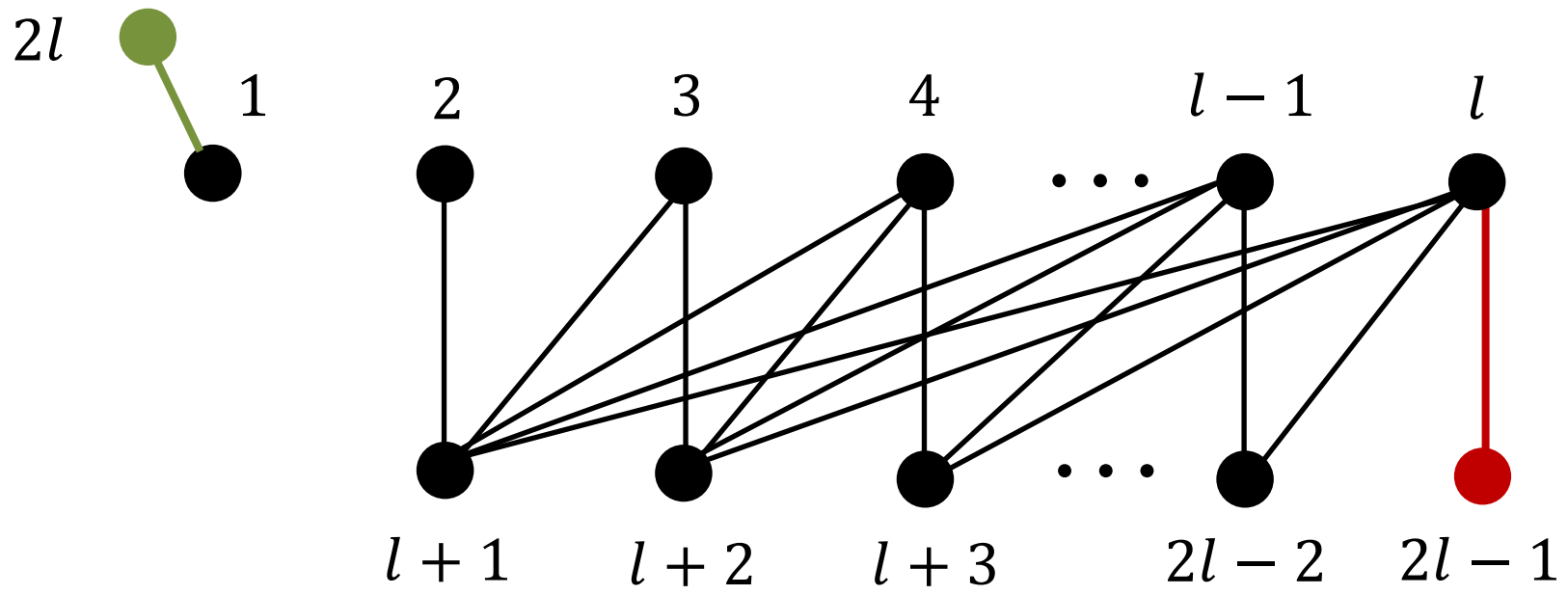

Therefore, to prove the inequality $x_{2l-1} - x_{2l} > 0$ it suffices to show that $|A| > |B|$. For this, we represent $A$ and $B$ as unions of disjoint sets:

$$A = A_1 \cup A_2 \text{ and } B = B_1 \cup B_2,$$

where $A_1 = \{H \in A | (l,k) \notin E(H) \ \forall k \in \{l+1, l+2, \dots, 2l-2\}\}$,

$A_2 = \{H \in A | \ \exists k \in \{l+1, l+2, \dots, 2l-2\}: (l,k) \in E(H)\}$,

$B_1 = \{G \in B | \ l \notin V(G)\}$, $B_2 = \{G \in B | \ l \in V(G)\}$,

and we shall prove that $|A_1| \geq |B_1|$, $|A_2| > |B_2|$ .

**1) We prove that $|A_1| \geq |B_1|$.**

Consider the map from the set $B_1$ to the set $A_1$, that assigns to each subgraph $G \in B_1$ with vertex set $V(G)$ and edge set $E(G)$ a subgraph from $A_1$ by the following rule:

$$2l \to 2l-1, \quad 1 \to l, \quad i \to i \text{ for } i \in V(G) \backslash \{1, 2l\},$$

$$(1, 2l) \to (l, 2l-1), \quad (1, i) \in E(G) \to (l, i) \text{ for } i \neq 2l,$$

$$(i, j) \in E(G) \to (i, j) \text{ for } i, j \neq 1.$$

Note that since vertex 1 in the graph $G$ can only be adjacent to vertices from the set $\{2, 3, \dots, l-1, 2l\}$, then under the map indicated the corresponding vertex $l$ can only be adjacent to vertices from the set $\{2, 3, \dots, l-1, 2l-1\}$. Thus, this map indeed maps graphs from the set $B_1$ into the set $A_1$. Clearly, this map is injective. Consequently, $|A_1| \geq |B_1|$.

**2) We prove that $|A_2| > |B_2|$.**

*Let us bound $|B_2|$ from above.*

Every subgraph from $B_2$ contains the vertices $1, l, 2l$ and does not contain the vertex $2l-1$. The remaining $n-3$ vertices for such subgraphs can be chosen in $\binom{2l-4}{n-3}$ ways. And since to each set of $n$ vertices of the graph $X_{2l}$ there correspond at most $n!$ subgraphs isomorphic to $F$ on this vertex set, we have

$$|B_2| \leq n! \binom{2l-4}{n-3}.$$

*Let us bound $|A_2|$ from below.*

Every subgraph from $A_2$ contains the vertices $l, 2l-1$. Next we consider only those subgraphs from $A_2$ in which the vertex $l$ is adjacent to exactly one of the $\frac{l}{2}$ vertices of the set $\boldsymbol{D} = \left\{l+1, l+2, \dots, \frac{3l}{2}\right\}$ (call this vertex $k$), and the remaining $n-3$ vertices lie among the $\frac{l}{2} - 1$ vertices of the set $\boldsymbol{U} = \left\{\frac{l}{2}+1, \frac{l}{2}+2, \dots, l-1\right\}$. Let $a$ be the number of such subgraphs.

Note that any 2 vertices from the set $\boldsymbol{U} \cup \{k\} \cup \{l\}$ are adjacent $X_{2l}$. Therefore, by adding to the fixed vertices $k, l, 2l-1$ and fixed edges $(l,k), (l, 2l-1)$ an arbitrary set of $n-3$ vertices from the set $U$, along with the necessary set of edges, one can form at least one subgraph isomorphic to $F$, in which the vertex $l$ is adjacent to exactly one vertex $k \in D$, and the vertex $2l-1$ is pendant.

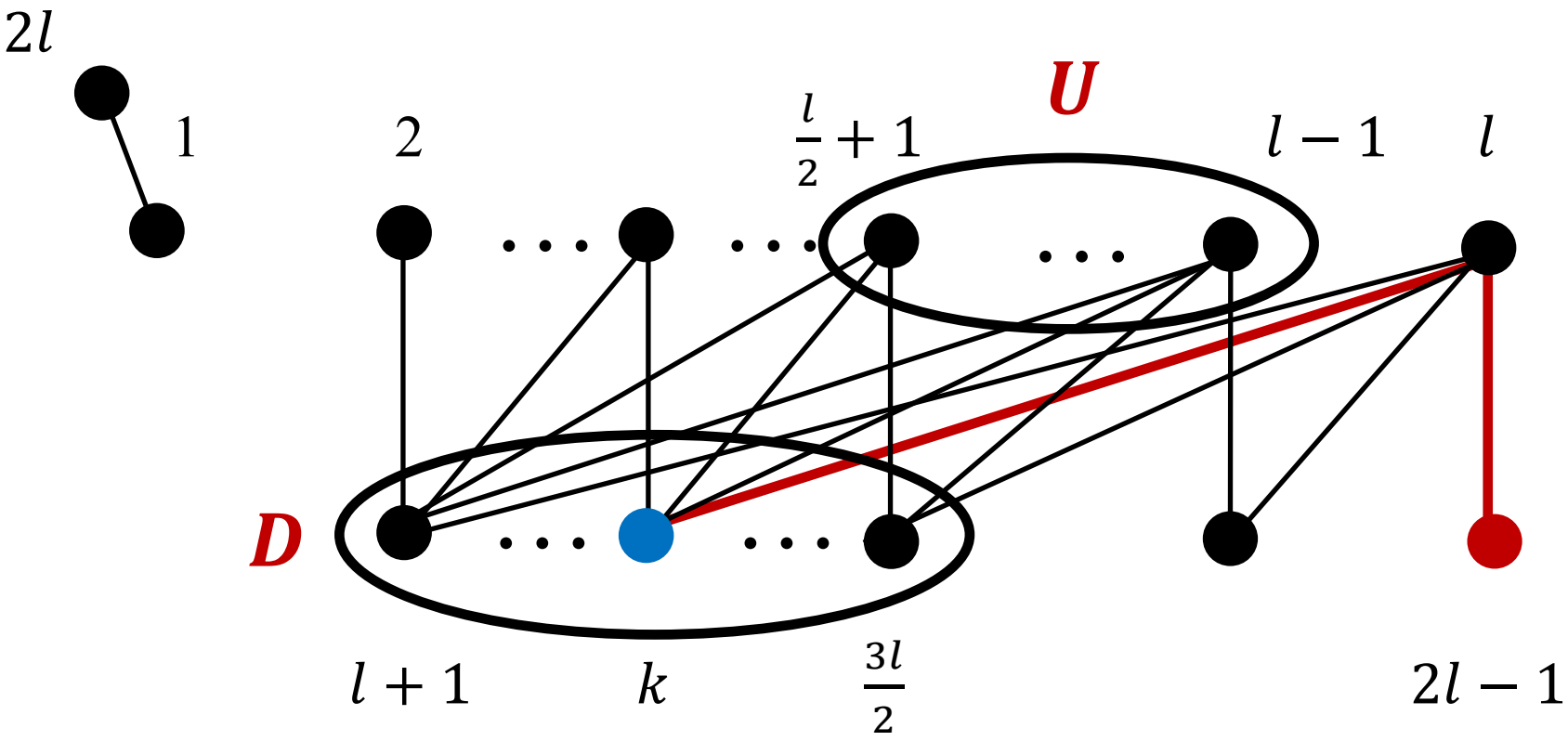

And since there are in total $\binom{\frac{l}{2}-1}{n-3}$ ways to choose $n-3$ vertices from the $\frac{l}{2}-1$ vertices of the set $\boldsymbol{U}$, we obtain at least $\binom{\frac{l}{2}-1}{n-3}$ subgraphs from $A_2$ for each $k$. It remains to note that $k$ can be chosen from the set $D$ in exactly $\frac{l}{2}$ ways. Therefore,

$$|A_2| \geq a \geq \frac{l}{2}\binom{\frac{l}{2}-1}{n-3}.$$

Taking into account the assumption $\frac{l}{2}\binom{\frac{l}{2}-1}{n-3} > n!\binom{2l-4}{n-3}$ and the inequality proved above $|B_2| \leq n!\binom{2l-4}{n-3}$, we obtain that $|A_2| > |B_2|$.

Then from 1)-2) it follows that $|A| > |B|$. □

**Proof of Lemma 2.6.** Let $l \in \mathbb{N}, l \vdots \operatorname{lcm}(m_1, m_2)$ and $l \geq \frac{m_1}{k_1}(x-b_1), l \geq \frac{m_2}{k_2}(y-b_2)$. Then the quantities $\binom{\frac{k_1}{m_1}l+b_1}{x}$ and $\binom{\frac{k_2}{m_2}l+b_2}{y}$ are defined. We consider 2 cases.

1) $y \in \mathbb{N}$. In this case, we have:

$$r_1\binom{\frac{k_1}{m_1}l+b_1}{x} > r_2\binom{\frac{k_2}{m_2}l+b_2}{y} \quad \Leftrightarrow$$

$$r_1 \cdot \frac{\left(\frac{k_1}{m_1}l+b_1\right)!}{x!\left(\frac{k_1}{m_1}l+b_1-x\right)!} - r_2 \cdot \frac{\left(\frac{k_2}{m_2}l+b_2\right)!}{y!\left(\frac{k_2}{m_2}l+b_2-y\right)!} > 0 \quad \Leftrightarrow$$

$$r_1 y!\left(\frac{k_1}{m_1}l+b_1-x+1\right)\left(\frac{k_1}{m_1}l+b_1-x+2\right)\cdot\ldots\cdot\left(\frac{k_1}{m_1}l+b_1\right) -$$

$$-r_2 x!\left(\frac{k_2}{m_2}l+b_2-y+1\right)\left(\frac{k_2}{m_2}l+b_2-y+2\right)\cdot\ldots\cdot\left(\frac{k_2}{m_2}l+b_2\right) > 0.$$

Next, consider the polynomial:

$$Q(t) = r_1 y! \left(\frac{k_1}{m_1}t + b_1 - x + 1\right) \cdot \left(\frac{k_1}{m_1}t + b_1 - x + 2\right) \cdot \ldots \cdot \left(\frac{k_1}{m_1}t + b_1\right) -$$

$$-r_2 x! \left(\frac{k_2}{m_2}t + b_2 - y + 1\right) \cdot \left(\frac{k_2}{m_2}t + b_2 - y + 2\right) \cdot \ldots \cdot \left(\frac{k_2}{m_2}t + b_2\right), t \in \mathbb{R}.$$

Note that $r_1 y! \left(\frac{k_1}{m_1}t + b_1 - x + 1\right) \cdot \left(\frac{k_1}{m_1}t + b_1 - x + 2\right) \cdot \ldots \cdot \left(\frac{k_1}{m_1}t + b_1\right)$ is a

polynomial of degree $x$ in the variable $t$, and

$r_2 x! \left(\frac{k_2}{m_2}t + b_2 - y + 1\right) \cdot \left(\frac{k_2}{m_2}t + b_2 - y + 2\right) \cdot \ldots \cdot \left(\frac{k_2}{m_2}t + b_2\right)$ is a polynomial

of degree $y$ in the variable $t$. Since $x > y$, it follows that $Q(t)$ is a polynomial of degree $x$ in the argument $t$. Moreover, the leading coefficient of the polynomial $Q(t)$, equal to $r_1 y! \cdot \frac{{k_1}^x}{{m_1}^x}$, is a positive number. Therefore, by a property of polynomials, there exists a number $M_1$ such that for all $t \geq M_1$ the polynomial $Q(t)$ takes only positive values.

In particular, $Q(l) > 0$ for all natural numbers $l \geq M_1,\ l \mathrel{\vdots} \operatorname{lcm}(m_1, m_2)$.

Let $N_1 = max\left(M_1, \frac{m_1}{k_1}(x - b_1), \frac{m_2}{k_2}(y - b_2)\right)$. Then, taking the above into account, for all natural numbers $l \geq N_1,\ l \mathrel{\vdots} \operatorname{lcm}(m_1, m_2)$ the following inequality holds:

$$r_1 \binom{\frac{k_1}{m_1}l + b_1}{x} > r_2 \binom{\frac{k_2}{m_2}l + b_2}{y}.$$

**2) $y = 0$.** In this case we have:

$$r_1 \binom{\frac{k_1}{m_1}l + b_1}{x} > r_2 \binom{\frac{k_2}{m_2}l + b_2}{0} \quad \Leftrightarrow$$

$$r_1 \cdot \frac{\left(\frac{k_1}{m_1}l + b_1\right)!}{x!\left(\frac{k_1}{m_1}l + b_1 - x\right)!} - r_2 > 0 \quad \Leftrightarrow$$

$$r_1 \left(\frac{k_1}{m_1}l + b_1 - x + 1\right) \cdot \left(\frac{k_1}{m_1}l + b_1 - x + 2\right) \cdot \ldots \cdot \left(\frac{k_1}{m_1}l + b_1\right) - r_2 x! > 0.$$

The further analysis of this inequality is analogous to Case 1. □

**Proof of Lemma 2.7.** Let $l \in \mathbb{N}, l \mathrel{\vdots} \operatorname{lcm}(m_1, m_2)$ and $l \geq \frac{m_1}{k_1}(x - b_1), l \geq \frac{m_2}{k_2}(x - b_2)$.

Then the quantities $\binom{\frac{k_1}{m_1}l + b_1}{x}$ and $\binom{\frac{k_2}{m_2}l + b_2}{x}$ are defined. Next, we consider 2 cases.

**1) $x = 0$.** In this case we have:

$$(p_1 l + r_1) \binom{\frac{k_1}{m_1}l + b_1}{0} > r_2 \binom{\frac{k_2}{m_2}l + b_2}{0} \quad \Leftrightarrow \quad p_1 l + r_1 > r_2.$$

Therefore, for all natural numbers $l \geq N_2 = max\left(\frac{r_2 - r_1}{p_1} + 1, -b_1\frac{m_1}{k_1}, -b_2\frac{m_2}{k_2}\right)$, $l \vdots \operatorname{lcm}(m_1, m_2)$ the following inequality holds:

$$(p_1 l + r_1)\binom{\frac{k_1}{m_1}l + b_1}{0} > r_2\binom{\frac{k_2}{m_2}l + b_2}{0}.$$

**2) $\boldsymbol{x} \in \mathbb{N}$.** In this case we have

$$(p_1 l + r_1)\binom{\frac{k_1}{m_1}l + b_1}{x} > r_2\binom{\frac{k_2}{m_2}l + b_2}{x} \quad \Leftrightarrow$$

$$(p_1 l + r_1) \cdot \frac{\left(\frac{k_1}{m_1}l + b_1\right)!}{x!\left(\frac{k_1}{m_1}l + b_1 - x\right)!} - r_2 \cdot \frac{\left(\frac{k_2}{m_2}l + b_2\right)!}{x!\left(\frac{k_2}{m_2}l + b_2 - x\right)!} > 0 \;\Leftrightarrow$$

$$(p_1 l + r_1)\left(\frac{k_1}{m_1}l + b_1 - x + 1\right) \cdot \left(\frac{k_1}{m_1}l + b_1 - x + 2\right) \cdot \ldots \cdot \left(\frac{k_1}{m_1}l + b_1\right) -$$

$$-r_2\left(\frac{k_2}{m_2}l + b_2 - x + 1\right) \cdot \left(\frac{k_2}{m_2}l + b_2 - x + 2\right) \cdot \ldots \cdot \left(\frac{k_2}{m_2}l + b_2\right) > 0.$$

Next, consider the polynomial:

$$P(t) = (p_1 t + r_1)\left(\frac{k_1}{m_1}t + b_1 - x + 1\right) \cdot \left(\frac{k_1}{m_1}t + b_1 - x + 2\right) \cdot \ldots \cdot \left(\frac{k_1}{m_1}t + b_1\right) -$$

$$-r_2\left(\frac{k_2}{m_2}t + b_2 - x + 1\right) \cdot \left(\frac{k_2}{m_2}t + b_2 - x + 2\right) \cdot \ldots \cdot \left(\frac{k_2}{m_2}t + b_2\right), \quad t \in \mathbb{R}.$$

Note that

$$(p_1 t + r_1)\left(\frac{k_1}{m_1}t + b_1 - x + 1\right) \cdot \left(\frac{k_1}{m_1}t + b_1 - x + 2\right) \cdot \ldots \cdot \left(\frac{k_1}{m_1}t + b_1\right)$$

is a polynomial of degree $x + 1$ in the variable $t$, and

$$r_2\left(\frac{k_2}{m_2}t + b_2 - x + 1\right) \cdot \left(\frac{k_2}{m_2}t + b_2 - x + 2\right) \cdot \ldots \cdot \left(\frac{k_2}{m_2}t + b_2\right)$$

is a polynomial of degree $x$ in the variable $t$. Therefore, $P(t)$ is a polynomial of degree $x + 1$ in $t$. Moreover, the leading coefficient of the polynomial $P(t)$, equal to $p_1 \cdot \frac{{k_1}^x}{{m_1}^x}$, is a positive number. Hence, by a property of polynomials, there exists a number $M_2$ such that for all $t \geq M_2$ the polynomial $P(t)$ takes only positive values. In particular, $P(l) > 0$ for all natural numbers $l \geq M_2$, $l \vdots lcm(m_1, m_2)$. Let

$$N_2 = max\left(M_2, \frac{m_1}{k_1}(x - b_1), \frac{m_2}{k_2}(x - b_2)\right).$$

Then, taking the above into account, for all natural numbers $l \geq N_2$, $l \vdots \operatorname{lcm}(m_1, m_2)$ the following inequality holds:

$$(p_1 l + r_1)\binom{\frac{k_1}{m_1}l + b_1}{x} > r_2\binom{\frac{k_2}{m_2}l + b_2}{x}. \ \square$$

**Proof of Lemma 3.1.** Consider a vertex $u \in V(F)$ with $\deg(u) = t$. If Lemma 3.1 does not hold, then in the graph $F$, taking into account its connectedness, there are no vertices other than $u$ and those adjacent to it. Thus, $n = t + 1$. But since $t$ is the minimum of the vertex degrees of $F$, all vertices in $F$ have degree $t = n - 1$, that is, $F$ is a complete graph. The diameter of a complete graph equals 1 ($n \geq 3$), while the diameter of $F$ is at least 2. Contradiction. Hence Lemma 3.1 holds. □

**Proof of Lemma 3.2**. Note that since the maximum label of a red vertex in $G_{2l-1}$ equals

$$1 + \frac{l}{n^2}\left(\frac{t(t-1)}{2} - 1\right)$$

and $t < n$, the label of any red vertex does not exceed

$$1 + \frac{l}{n^2}\left(\frac{t(t-1)}{2} - 1\right) \; < \; 1 + \frac{l}{n^2} \cdot \frac{t^2}{2} \; < \; 1 + \frac{l}{n^2} \cdot \frac{n^2}{2} \; = \; 1 + \frac{l}{2}.$$

Consequently, the label of a red vertex is at most $\frac{l}{2}$. □

**Proof of Lemma 3.3.** Let $l \vdots 2n^2$ and $\binom{\frac{l}{2}-1}{n-2} > n!\,\binom{2l-t-3}{n-t-1}$. Consider in the graph $G_{2l-1}$ two consecutive vertices $i$ and $i + 1$, where $i \in \{1, 2, \dots, l - 1\}$. We estimate the difference $g_{i+1} - g_i$. Note that the vertices $i$ and $i + 1$ cannot both be red, since otherwise $\frac{l}{n^2} = 1$, whence $l = n^2$, which contradicts the condition $l \vdots 2n^2$. Hence, three cases are possible.

**1) The vertices $i$ and $i + 1$ are black.**

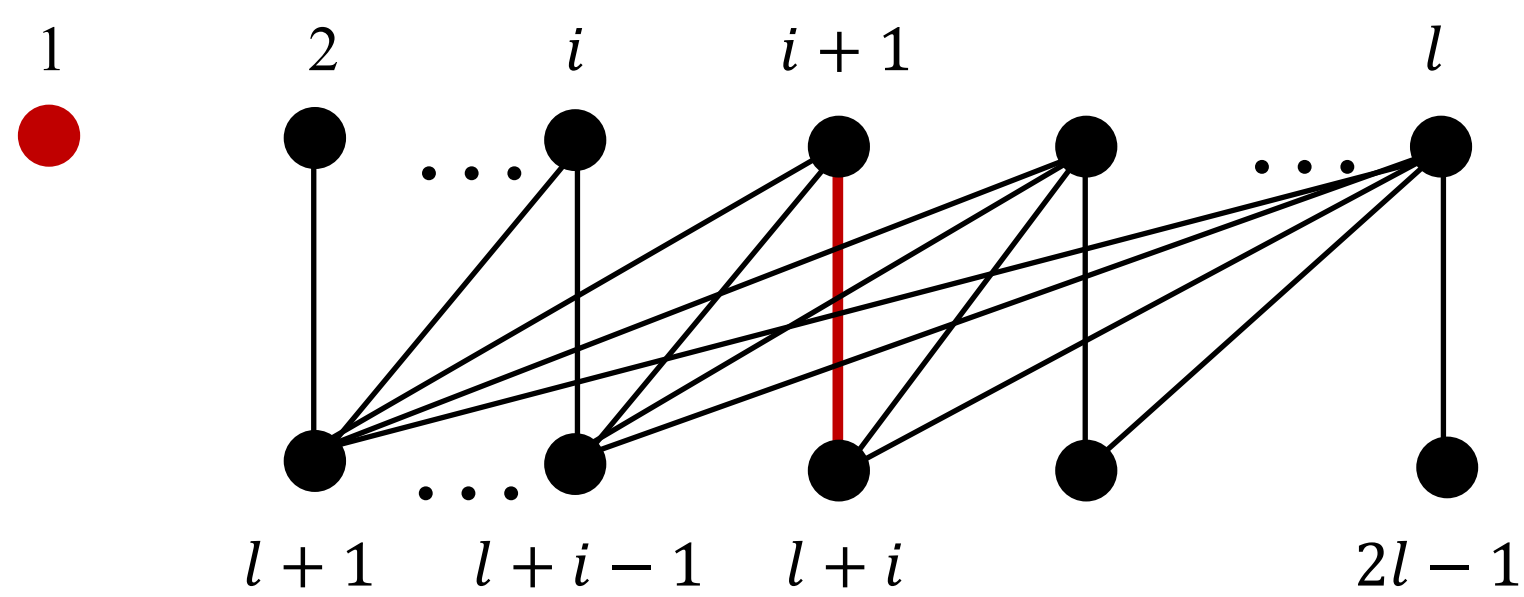

Note that if we delete the edge $(i + 1, l + i)$ from $G_{2l-1}$, then the vertices $i$ and $i + 1$ become symmetric, since both are adjacent to the vertices $l + 1, l + 2, \dots, l + i - 1$, as well as to all vertices of the upper level other than $i$ and $i + 1$. Moreover, the vertices $i$ and $i + 1$ in the graph $G_{2l-1} \setminus \{(i + 1, l + i)\}$ are not adjacent to the vertices $l + i, l + i + 1, \dots, 2l - 1$. Consequently, the $F$-degrees of the vertices $i + 1$ and $i$ in $G_{2l-1} \setminus \{(i + 1, l + i)\}$ coincide. Hence the difference of the $F$-degrees $(g_{i+1} - g_i)$ of these vertices in $G_{2l-1}$ equals the number of subgraphs isomorphic to $F$ that contain the edge $(i + 1, l + i)$ and do not contain the vertex $i$. The set of such subgraphs is nonempty. Let us prove this.

Note that the degree of the vertex $l + i$ is at least $t$. Consider the vertex $l + i$, together with $t$ vertices of the upper level of $G_{2l-1}$ adjacent to it, including the vertex $i + 1$. From the definition of $G_{2l-1}$ and from the fact that the vertex $i$ is black, it follows that the vertex $i$ is not among these $t + 1$ vertices. Next, add to these vertices $n - t - 1$ arbitrary new vertices of the upper level, other than $i$. From the resulting set of vertices and the edge $(i + 1, l + i)$, by adding the necessary edges one can form at least one subgraph isomorphic to $F$ and containing the edge $(i + 1, l + i)$, in which the vertex $l + i$ has minimum degree $t$. By construction, such a subgraph of $G_{2l-1}$ does not contain the vertex $i$. Hence, $g_{i+1} - g_i > 0$.

**2) The vertex $i$ is black, the vertex $i+1$ is red.**

Let the red vertex $i+1$ be adjacent to some vertex $j \in \{2l-t+1, 2l-t+2, \dots, 2l-1\}$. Similarly to case 1) it is shown that the $F$-degrees of the vertices $i+1$ and $i$ in the graph $G_{2l-1} \setminus \{(i+1, l+i), (i+1, j)\}$ coincide. Hence the difference of the $F$-degrees $(g_{i+1} - g_i)$ of these vertices in $G_{2l-1}$ equals the number of subgraphs isomorphic to $F$ that contain the edge $(i+1, l+i)$ or the edge $(i+1, j)$ and do not contain the vertex $i$.

The set of such subgraphs is nonempty, which is established in the same way as in case 1).

**3) The vertex $i$ is red, the vertex $i+1$ is black.**

Let the red vertex $i$ be adjacent to some vertex $j \in \{2l-t+1, 2l-t+2, \dots, 2l-1\}$. As in case 1), it is shown that the $F$-degrees of the vertices $i+1$ and $i$ in the graph $G_{2l-1} \setminus \{(i+1, l+i), (i, j)\}$ coincide.

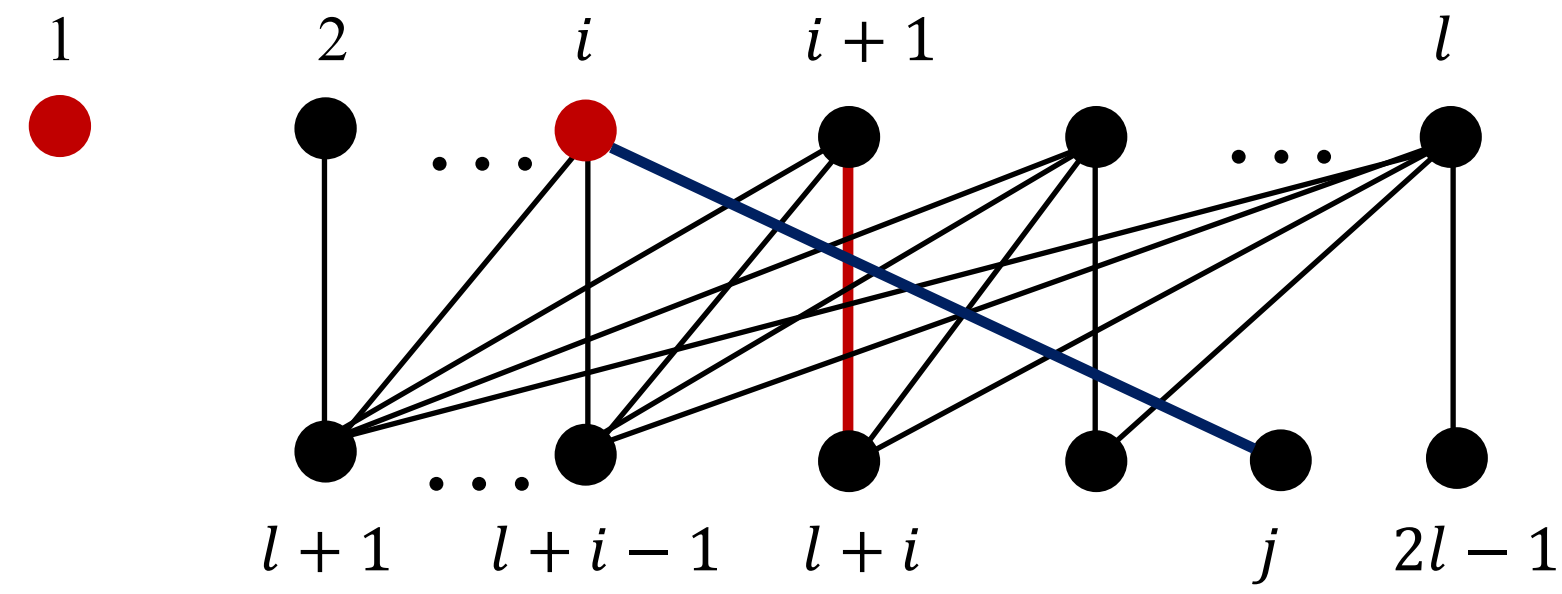

Then $g_{i+1} - g_i = |X| - |Y|$, where $X$ is the set of subgraphs isomorphic to $F$ that contain the edge $(i+1, l+i)$ and do not contain the vertex $i$, and $Y$ is the set of subgraphs isomorphic to $F$ that contain the edge $(i, j)$ and do not contain the vertex $i+1$.

*We estimate* $|Y|$ *from above.*

Since $l \in N$ and $l \vdots 2n^2$, we have $l \geq 2n^2$. Then $t < n < n^2 \leq \frac{l}{2}$. Moreover, from Lemma 3.2 it follows that $i \leq \frac{l}{2}$. Therefore,

$$j \geq 2l - t + 1 > 2l - \frac{l}{2} + 1 > l + i.$$

In particular, from the last inequality it follows that $j$ is not adjacent to $i+1$.

Further, since the vertex $j$ has degree $t$ in $G_{2l-1}$, every subgraph in $Y$ contains, together with the vertex $j$, all $t$ vertices adjacent to $j$. Moreover, such subgraphs do not contain the vertex $i+1$. Hence the remaining $n-t-1$ vertices for such subgraphs can be chosen in $\binom{2l-t-3}{n-t-1}$ ways. Since to each set of $n$ vertices of the graph $G_{2l-1}$ there correspond at most $n!$ subgraphs isomorphic to $F$ on this vertex set, it follows that

$$|Y| \leq n! \binom{2l-t-3}{n-t-1}.$$

*We estimate* $|X|$ *from below.*

Note that since $i$ is a red vertex, by Lemma 3.2 the inequality $i \leq \frac{l}{2}$ holds.
Then $l - (i+1) \geq \frac{l}{2} - 1$.

Every subgraph in $X$ contains the edge $(i+1, l+i)$ and does not contain the vertex $i$. Next, consider only those subgraphs in $X$ for which the remaining $n-2$ vertices lie among the $l-(i+1)$ vertices of the set $\{i+2, i+3, \dots, l\}$. Let $x$ be the number of such subgraphs. Note that since any two vertices of the set $\{i+1, i+2, i+3, \dots, l, l+i\}$ are adjacent in $G_{2l-1}$, by adding to the vertices $i+1$, $l+i$ and the edge $(i+1, l+i)$ an arbitrary set of $n-2$ vertices of the set $\{i+2, i+3, \dots, l\}$, together with the necessary set of edges, one can form at least one subgraph isomorphic to $F$. Since there are in total $\binom{l-(i+1)}{n-2}$ ways to choose $n-2$ vertices from the $l-(i+1)$ vertices, it follows that

$$x \geq \binom{l-(i+1)}{n-2} \geq \binom{\frac{l}{2}-1}{n-2}.$$

Consequently,

$$|X| \geq x \geq \binom{\frac{l}{2}-1}{n-2}.$$

And, taking into account the assumption $\binom{\frac{l}{2}-1}{n-2} > n! \binom{2l-t-3}{n-t-1}$ and the inequality $|Y| \leq n! \binom{2l-t-3}{n-t-1}$ proved above, we obtain

$$g_{i+1} - g_i = |X| - |Y| \geq \binom{\frac{l}{2}-1}{n-2} - n! \binom{2l-t-3}{n-t-1} > 0. \ \square$$

**Proof of Lemma 3.4**. Let $l \vdots 2n^2$. Consider in the graph $G_{2l-1}$ two vertices 1 and $l+1$. We estimate the difference $g_1 - g_{l+1}$. As in case 1) of the proof of Lemma 2.3, it is established that the $F$-degrees of the vertices 1 and $l+1$ in the graph $G_{2l-1} \setminus \{(1, 2l-1)\}$ coincide. Hence the difference of the $F$-degrees $(g_1 - g_{l+1})$ of these vertices in $G_{2l-1}$ equals the number of subgraphs isomorphic to $F$ that contain the edge $(1, 2l-1)$ and do not contain the vertex $l+1$. The set of such subgraphs is nonempty. Let us prove this.

Consider the vertex $2l-1$, together with all $t$ vertices of the upper level of the graph $G_{2l-1}$ adjacent to it, among which, evidently, is the vertex 1. Moreover, add to these vertices $n-t-1$ arbitrary new vertices of the upper level. From the resulting set of vertices and the edge $(1, 2l-1)$, by adding the necessary edges one can always form at least one subgraph isomorphic to $F$ and containing the edge $(1, 2l-1)$, in which the vertex $2l-1$ has minimum degree $t$. By construction, such a subgraph of $G_{2l-1}$ does not contain the vertex $l+1$. Consequently, $g_1 - g_{l+1} > 0$. □

**Proof of Lemma 3.5**. Let $l \vdots 2n^2$. Consider in the graph $G_{2l-1}$ two consecutive vertices $i$ and $i+1$, where $i \in \{l+1, l+2, \dots, 2l-t-1\}$.

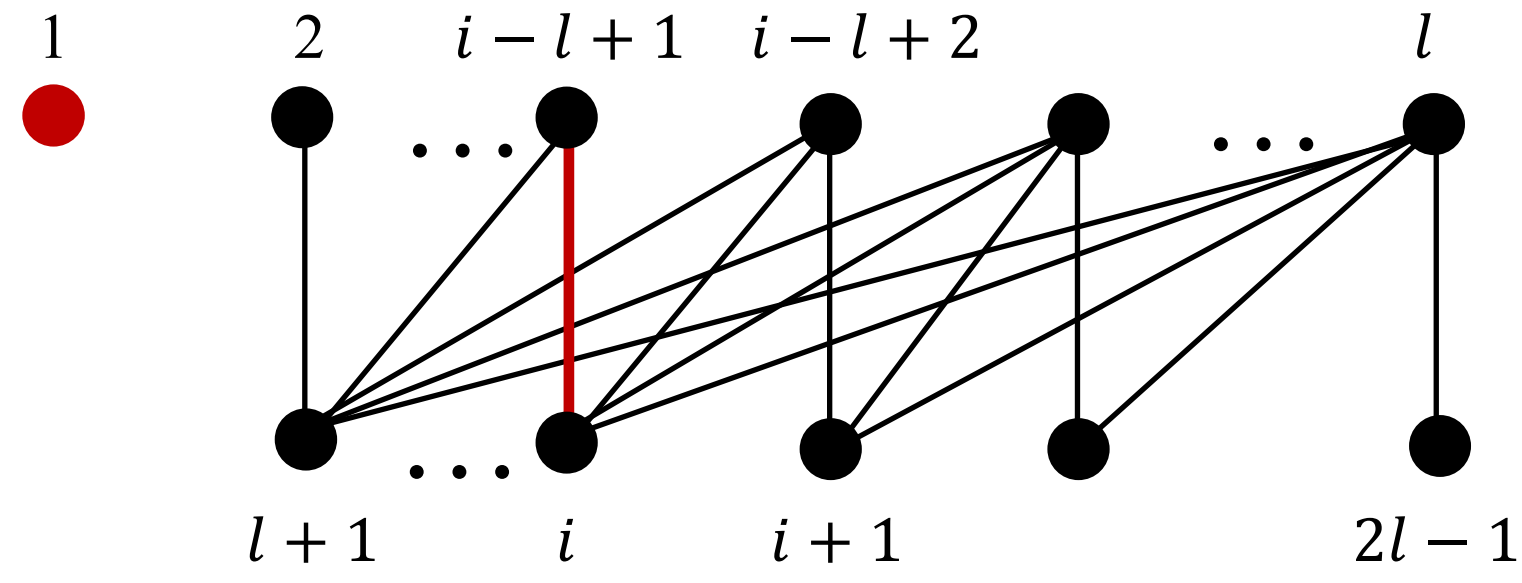

We estimate the difference $g_i - g_{i+1}$. If we delete the edge $(i, i-l+1)$ from $G_{2l-1}$ then the vertices $i$ and $i+1$ become symmetric, since both are adjacent to the vertices $i-l+2$, $i-l+3$, …, $l$, and are adjacent to no other vertices. Consequently, the $F$-degrees of the vertices $i$ and $i+1$ in the graph $G_{2l-1} \setminus \{(i, i-l+1)\}$ coincide. Hence the difference of the $F$-degrees $(g_i - g_{i+1})$ of these vertices in $G_{2l-1}$ equals the number of subgraphs isomorphic to $F$ that contain the edge $(i, i-l+1)$ and do not contain the vertex $i+1$. The set of such subgraphs is nonempty. Let us prove this.

Note that the degree of the vertex $i$ is greater than $t$. Consider the vertex $i$, together with $t$ vertices of the upper level of $G_{2l-1}$ adjacent to it, including the vertex $i-l+1$. Moreover, add to these vertices $n-t-1$ arbitrary new vertices of the upper level. From the resulting set of vertices and the edge $(i, i-l+1)$, by adding the necessary edges one can form at least one subgraph isomorphic to $F$ and containing the edge $(i, i-l+1)$, in which the vertex $i$ has minimum degree $t$. By construction, such a subgraph of $G_{2l-1}$ does not contain the vertex $i+1$. Consequently, $g_i - g_{i+1} > 0$. □

**Proof of Lemma 3.6.** Let $l \vdots 2n^2$ and $\frac{l}{n^2}\binom{\frac{l}{n^2}}{n-t-2} > n!\,(t-1)\binom{2l-t-4}{n-t-2}$. Consider in the graph $G_{2l-1}$ two consecutive vertices $i$ and $i+1$, where $i \in \{2l-t, 2l-t+1, \dots, 2l-2\}$. According to the structure of $G_{2l-1}$, the vertex $i$ is adjacent in this graph to exactly $2l-i$ black vertices $i-l+1, i-l+2, \dots, l$, and to exactly $t+i-2l$ red vertices.

Note that the vertex $i = 2l-t$ is not joined by edges to any red vertex.

In the case $i \geq 2l-t+1$, let $m_1, m_2, \dots, m_{t+i-2l}$ be the labels of the red neighbors of the vertex $i$ in increasing order.

Let $i \geq 2l-t$. The vertex $i+1$ is adjacent to exactly $2l-i-1$ black vertices, namely $i-l+2, i-l+3, \dots, l$, and to exactly $t+i-2l+1$ red vertices with labels $k_1, k_2, \dots, k_{t+i-2l+1}$ in increasing order.

Let $i \geq 2l-t+1$. Then, by the definition of $G_{2l-1}$, the following inequalities hold:

$$k_1 < k_2 < \cdots < k_{t+i-2l+1} < m_1 < m_2 < \cdots < m_{t+i-2l},$$

with

$$m_1 - k_{t+i-2l+1} = \frac{l}{n^2}.$$

To simplify notation, denote $k_{t+i-2l+1} = a$, $m_1 = b$. Then $b - a = \frac{l}{n^2}$.

Further, in the case $i = 2l-t$ as well, we set $a = k_1$, $b = a + \frac{l}{n^2}$, but in this case the vertex $b$ will be black.

**Assertion 3.1.**

*Let* $, n, t \in N, n \geq 3, l \vdots 2n^2$, $t < n$. *Then the following inequalities hold:*

1) $t < \dfrac{l}{6}$;

2) $l - b - t > \dfrac{l}{n^2}$ $\;for\;$ $i \geq 2l-t$;

3) $b < i - l + 1$ $\;for\;$ $i \geq 2l-t$;

4) $m_{t+i-2l} < i - l + 1$ $\;for\;$ $i \geq 2l-t+1$.

**Proof of 1):** Since $\in N, l \vdots 2n^2, t < n, n \geq 3$, we have

$$l \geq 2n^2 > 2nt \geq 6t,$$

whence $t < \frac{l}{6}$.

**Proof of 2):** Since $a \leq \frac{l}{2}$ (by Lemma 3.2), $t < \frac{l}{6}$ (by Assertion 3.1.(1)), $n \geq 3$, we have

$$l - b - t = l - a - \frac{l}{n^2} - t > l - \frac{l}{2} - \frac{l}{9} - \frac{l}{6} = \frac{2}{9}l > \frac{l}{n^2}.$$

**Proof of 3):** From Assertion 3.1.(2) it follows that $b < l - t - \frac{l}{n^2}$. But, since $i \geq 2l - t$, we have

$$l - t - \frac{l}{n^2} < l - t + 1 \leq i - l + 1.$$

Consequently, $b < i - l + 1$ for $i \geq 2l - t$.

**Proof of 4):** By Lemma 3.2, $m_{t+i-2l} \leq \frac{l}{2}$. Moreover, by Assertion 3.1.(1), we have $t < \frac{l}{6}$. Then for $i \geq 2l - t + 1$ the following inequality holds:

$$i - l + 1 \geq 2l - t - l + 2 > l - t > l - \frac{l}{6} > \frac{l}{2} \geq m_{t+i-2l}.$$

Hence, $m_{t+i-2l} < i - l + 1$ for $i \geq 2l - t + 1$. Assertion 3.1 is proved.

From Assertions 3.1.(3) and 3.1.(4) it follows that for $i \geq 2l - t$, the vertex $b$ and all red vertices adjacent to $i$ or $i + 1$ lie on the upper level of the graph $G_{2l-1}$ to the left of the vertex $i - l + 1$.

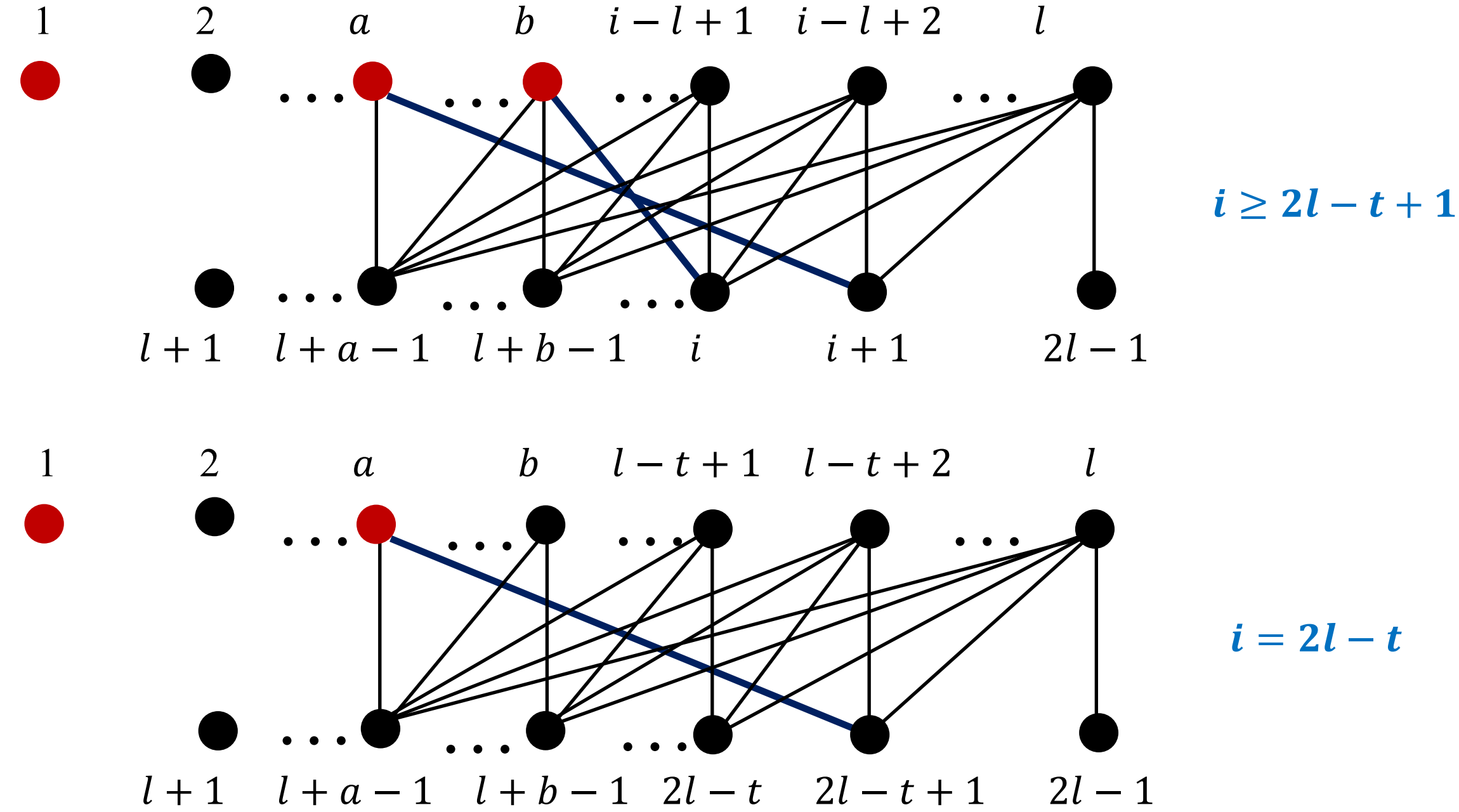

We prove that $g_i - g_{i+1} > 0$. Clearly, the quantity $g_i - g_{i+1}$ equals the difference between the number of subgraphs isomorphic to $F$ that contain the vertex $i$ and do not contain the vertex $i + 1$ (denote the set of such subgraphs by $S$), and the number of subgraphs isomorphic to $F$ that contain the vertex $i + 1$ and do not contain the vertex $i$ (denote the set of such subgraphs by $T$). Therefore, to prove $g_i - g_{i+1} > 0$ it suffices to show that $|S| > |T|$.

Note that since the vertex $i$ has degree $t$ in $G_{2l-1}$, in every subgraph from the set $S$ the vertex $i$ is adjacent to the black vertices $i-l+1, i-l+2, \dots, l$, and, in the case $i \geq 2l-t+1$, to the red vertices $m_1, m_2, \dots, m_{t+i-2l}$. Similarly, the vertex $i+1$ in every subgraph from the set $T$ is adjacent to the black vertices $i-l+2, i-l+3, \dots, l$, and to the red vertices $k_1, k_2, \dots, k_{t+i-2l+1}$.

Next, represent $S$ and $T$ as unions of disjoint sets:

$$S = S_1 \cup S_2 \quad and \quad T = T_1 \cup T_2,$$

where

$$S_1 = \{H \in S \mid (i-l+1, j) \notin E(H) \ \ \forall j \in \{l+a, l+a+1, \dots, l+b-1\}\},$$

$$S_2 = \{H \in S \mid \exists j \in \{l+a, l+a+1, \dots, l+b-1\} : \ (i-l+1, j) \in E(H)\},$$

$$T_1 = \{G \in T \mid V(G) \cap \{m_1, m_2, \dots, m_{t+i-2l}, i-l+1\} \ = \ \emptyset\},$$

$$T_2 = \{G \in T \mid V(G) \cap \{m_1, m_2, \dots, m_{t+i-2l}, i-l+1\} \neq \emptyset\},$$

where the set $\{m_1, m_2, \dots, m_{t+i-2l}\}$ is empty if $i = 2l-t$.

Let us proof that $|S_1| \geq |T_1|$ and $|S_2| > |T_2|$.

**1) We prove that $|S_1| \geq |T_1|$.**

Let $\{i+1, k_1, k_2, \dots, k_{t+i-2l+1}\} \ = \ K$. Consider a mapping of the set $T_1$ into the set $S_1$, under which to each subgraph $G \in T_1$ with vertex set $V(G)$ and edge set $E(G)$ we assign a subgraph in $S_1$ according to the following rule:

$$i+1 \to i, \ \ k_p \to m_p \ \ for \ \ p \in \{1, 2, \dots, t+i-2l\}, \ \ k_{t+i-2l+1} = a \to i-l+1,$$

$$s \to s \ \ for \ \ s \in V(G) \setminus K,$$

$$(i+1, k_p) \to (i, m_p) \ \ for \ \ p \in \{1, 2, \dots, t+i-2l\},$$

$$(i+1, k_{t+i-2l+1}) \to (i, i-l+1),$$

$$(i+1, s) \in E(G) \to (i, s) \ \ for \ \ s \in V(G) \setminus K,$$

$$(s, k_p) \in E(G) \to (s, m_p) \ \ for \ \ p \in \{1, 2, \dots, t+i-2l\}, s \in V(G) \setminus K,$$

$$(s, k_{t+i-2l+1}) \in E(G) \to (s, i-l+1) \ \ for \ \ s \in V(G) \setminus K,$$

$$(s, z) \in E(G) \to (s, z) \ \ for \ \ s, z \in V(G) \setminus K,$$

$$(k_q, k_p) \in E(G) \to (m_q, m_p) \ \ for \ \ p, q \in \{1, 2, \dots, t+i-2l\},$$

$$(k_p, k_{t+i-2l+1}) \in E(G) \to (m_p, i-l+1) \ \ for \ \ p \in \{1, 2, \dots, t+i-2l\}.$$

Note that since the vertex $k_{t+i-2l+1} = a$ in the graph $G$ is not adjacent to any vertex of the set $\{l+a, l+a+1, \dots, l+b-1\}$, under this mapping the corresponding vertex $i-l+1$ is not adjacent to any vertex of $\{l+a, l+a+1, \dots, l+b-1\}$. Hence, this mapping indeed maps graphs from the set $T_1$ into the set $S_1$.

Clearly, this mapping is injective. Consequently, $|S_1| \geq |T_1|$.

**2) We prove that $|S_2| > |T_2|$.**

*We estimate $|T_2|$ from above.*

Every subgraph in $T_2$ contains the vertex $i+1$ and $t$ fixed vertices of the upper level adjacent to $i+1$. Moreover, every such subgraph contains a vertex from the set $\{m_1, m_2, \dots, m_{t+i-2l}, i-l+1\}$ (there are at most $t-1$ choices for such a vertex, since $t+i-2l \le t-2$ for $i \le 2l-2$) and does not contain the vertex $i$.
The remaining $n-t-2$ vertices for such subgraphs can be chosen in $\binom{2l-t-4}{n-t-2}$ ways.

Since to each set of $n$ vertices of the graph $G_{2l-1}$ there correspond at most $n!$ subgraphs isomorphic to $F$ on this vertex set, it follows that

$$|T_2| \le n!\,(t-1)\binom{2l-t-4}{n-t-2}.$$

*We estimate $|S_2|$ from below.*

In every subgraph of $S_2$ the vertex $i$ is adjacent to exactly $t$ vertices of the upper level. Next, consider only those subgraphs of $S_2$ for which the vertex $i-l+1$ is adjacent to exactly one of the $\frac{l}{n^2}$ vertices of the set

$$D = \{l+a, l+a+1, \dots, l+b-1\}$$

(call this vertex $j$), and the remaining $n-t-2$ vertices lie among the vertices of the set

$$U = \{b, b+1, \dots, l\} \setminus N(i),$$

where $N(i)$ is the set of vertices adjacent to $i$. Let $q$ be the number of such subgraphs.

From Assertions 3.1.(3) and 3.1.(2) it follows that $N(i) \subseteq \{b, b+1, \dots, l\}$ and

$$|U| = l-b+1-t > \frac{l}{n^2}.$$

Note that since any two vertices of the set $\{b, b+1, \dots, l\} \cup \{j\}$ are adjacent in $G_{2l-1}$, by adding to the fixed vertices $i, N(i), j$ and the fixed edge $(i-l+1, j)$ an arbitrary set of $n-t-2$ vertices from the set $U$, together with the necessary set of edges, one can form at least one subgraph isomorphic to $F$ in which the vertex $i-l+1$ is adjacent to exactly one vertex $j \in D$, and the vertex $i$ has minimum degree $t$. Such a subgraph indeed exists, by Lemma 3.1, where the roles of the vertices $u, v, w$ are played by the vertices $i, i-l+1, j$, respectively.

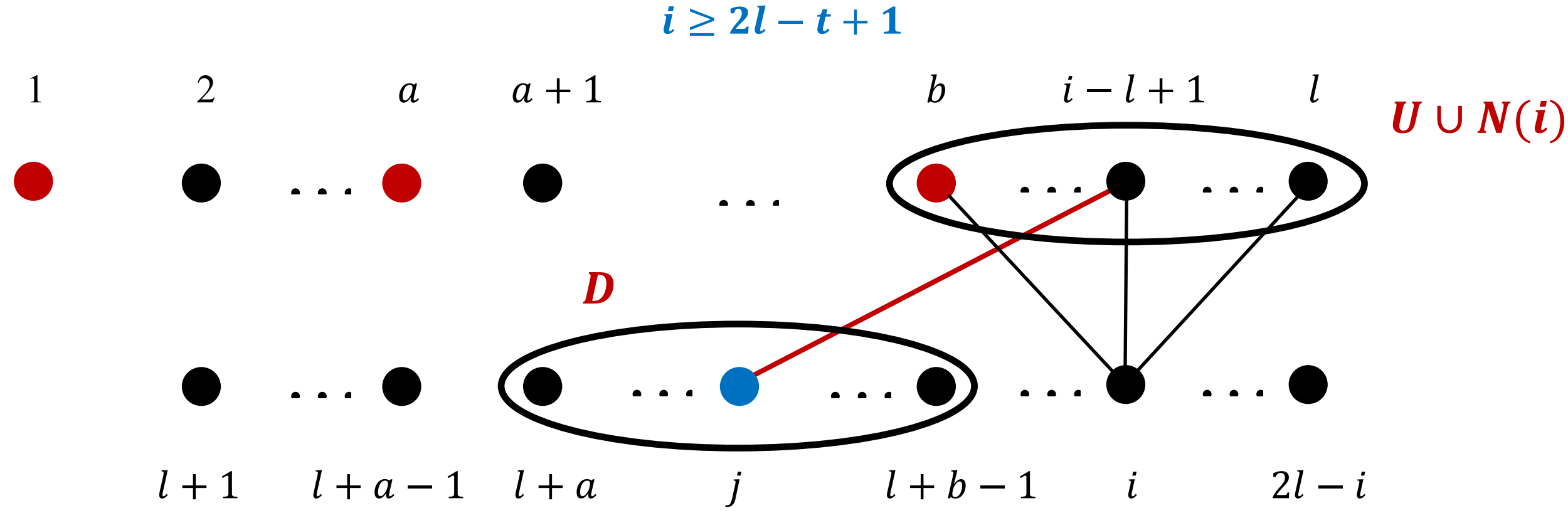

Since there are in total $\binom{|U|}{n-t-2} \geq \binom{\frac{l}{n^2}}{n-t-2}$ ways to choose $n-t-2$ vertices from the $|U|$ vertices of the set $U$, we obtain at least $\binom{\frac{l}{n^2}}{n-t-2}$ subgraphs from $S_2$ for each $j$. It remains to note that $j$ can be chosen from the set $D$ in exactly $\frac{l}{n^2}$ ways. Consequently,

$$|S_2| \geq q \geq \frac{l}{n^2}\binom{\frac{l}{n^2}}{n-t-2}.$$

And, taking into account the assumption

$$\frac{l}{n^2}\binom{\frac{l}{n^2}}{n-t-2} \; > \; n!\,(t-1)\binom{2l-t-4}{n-t-2}$$

and the inequality proved above $|T_2| \leq n!\,(t-1)\binom{2l-t-4}{n-t-2}$, we obtain $|S_2| > |T_2|$.

Finally, from 1)–2) it follows that $|S| > |T|$. □